\documentclass[11pt, a4paper]{amsart}

\usepackage{proof-at-the-end}

\usepackage[table]{xcolor}
\usepackage{amssymb, xargs}
\usepackage{amsthm}
\usepackage{amsaddr}
\usepackage{subcaption}
\usepackage{tikz}
\usepackage{amsmath}
\usepackage{mathtools}
\usepackage{chngcntr}
\usepackage[utf8]{inputenc}
\usepackage[sort, nocompress]{cite}
\usepackage[margin=0.75in]{geometry}
\usepackage{enumitem}

\usepackage{fancyhdr}
\usepackage{epigraph}
\usepackage{cleveref}

\newtheorem{theorem}{Theorem}[section]
\newtheorem*{result}{Main Result}

\newtheorem{assumption}[theorem]{Assumption}
\newtheorem{remark}[theorem]{Remark}

\newtheorem{example}[theorem]{Example}

\counterwithin{figure}{section}
\counterwithin{table}{section}
\counterwithin{equation}{section}
\counterwithin{theorem}{section}

\newcommand{\RefRem}[1]{}

\newcommand{\mc}{\mathcal}

\newcommand{\bs}{\boldsymbol}
\renewcommand{\H}{\bs{\mc{H}}}
\newcommand{\X}{\bs{\mc{H}}}
\newcommand{\R}{\ensuremath{\mathbb{R}}}
\newcommand{\C}{\ensuremath{\mathbb{C}}}
\newcommand{\N}{\ensuremath{\mathbb{N}}}

\newcommandx{\norm}[2][2=]{\ensuremath{\left\| #1 \right\|_{#2}}}

\newcommand{\dom}[1]{\mathcal{D}\mathopen{}\left(#1\right)\mathclose{}}
\renewcommand{\to}[1]{\mathcal{T}\mathopen{}\left(#1\right)\mathclose{}}
\renewcommand{\oe}[3]{\to{\exp{\int_{#2}^{#3}#1(\tau)\d\tau}}}

\renewcommand{\i}{\ensuremath{\mathrm{i}}}
\newcommand{\id}{\ensuremath{I}}
\newcommand{\inp}[2]{\mathopen{}\left\langle #1,\, #2\right\rangle\mathclose{}}
\renewcommand{\d}{\mathop{}\,\mathrm{d}}
\newcommand{\uu}{\ensuremath{\bs{\psi}}}
\newcommand{\Lnorm}{\ensuremath{\mc{L}}}
\newcommand{\comm}[2]{\mathopen{}\left[ #1,\, #2\right]\mathclose{}}
\newcommand{\Lp}[3][2]{\ensuremath{L^{#1}\mathopen{}\left(#2,\,#3\right)\mathclose{}}}
\newcommand{\set}[1]{\ensuremath{\mathopen{}\left\{#1\right\}\mathclose{}}}
\newcommand{\LL}[1]{\ensuremath{\mathcal{L}\mathopen{}\left(#1\right)\mathclose{}}}
\newcommand{\diff}{\Theta}
\renewcommand{\exp}[1]{\operatorname{exp}\mathopen{}\left[#1\right]\mathclose{}}
\newcommand{\adj}[1]{{#1}^*}
\DeclareMathOperator*{\dist}{dist}

\newcommand{\normsq}[2]{\ensuremath{\mathopen{}\left\|#1\right\|^2_{#2}\mathclose{}}}
\newcommand{\OB}{\ensuremath{\tilde{B}}}
\newcommand{\supp}[1]{\ensuremath{\operatorname{supp}\mathopen{}\left(#1\right)\mathclose{}}}

\title{Ranks of Tensor Networks for Eigenspace Projections and the Curse of
Dimensionality}
\author{Mazen Ali}
\thanks{Parts of this work were completed while the author was at
the Institute for Numerical Mathematics, Ulm University, Germany.}
\address{Centrale Nantes, LMJL UMR CNRS 6629, France}
\email{mazen.ali@ec-nantes.fr}


\begin{document}

\begin{abstract}
    The hierarchical (multi-linear) rank of
    an order-$d$ tensor
    is key in determining the cost of representing a tensor as a
    (tree) Tensor Network (TN).
    In general, it is known that,
    for a fixed accuracy,
    a tensor with random entries
    cannot be expected to be efficiently
    approximable without the \emph{curse of dimensionality},
    i.e., a complexity growing exponentially with $d$.
    In this work, we show that the ground state projection (GSP)
    of a class of unbounded Hamiltonians can be
    approximately represented as an operator of low
    effective dimensionality that is independent of the (high) dimension $d$
    of the GSP.
    This allows to approximate the GSP
    \emph{without} the curse of dimensionality.
\end{abstract}

\maketitle

\section{Introduction}
Tensor product methods naturally arise in high-dimensional problems, e.g.,
when one is interested in representing $d$-dimensional
arrays $\bs\psi\in\R^{n\times\ldots\times n}$ or approximating
$d$-dimensional functions
$\bs\psi:\R^d\rightarrow\R$,
and their application is ubiquitous across many scientific fields.
We refer to \cite{HB,TN,Khoromskij18,Khoromskaia} for overviews
and examples.

In the case of an array $\bs\psi\in\R^{n^d}$, the number of entries
$n^d$ scales exponentially with $d$
-- an expression of a phenomenon commonly
referred to as the
\emph{curse of dimensionality}\footnote{Coined by R. Bellman \cite{Bellman:DynamicProgramming}.}.
In order to recover \emph{tractable}\footnote{For a definition
of tractability see \cite{WOZNIAKOWSKI199496}.} computational methods,
a common technique is tensor product approximation
\begin{align}\label{eq:cp}
    \bs\psi\approx\bs\psi_r=\sum_{k=1}^rv^1_k\otimes\ldots\otimes v_k^d,
\end{align}
where $r$ is referred to as the \emph{tensor rank} of
$\bs\psi_r$.
Here, the representation complexity of $\bs\psi_r$ is
$ndr$, i.e., linear in $d$ and $r$.
The utility of such an approximation, of course,
depends on how $r$ scales w.r.t.\ $d$ for a given
approximation accuracy $\varepsilon>0$
\begin{align}\label{eq:basicq}
    r_\varepsilon:=
    \min
    \left\{r\in\N_0:\;\norm{\bs\psi-\bs\psi_r}\leq\varepsilon\right\}=
    \mc O(?),
\end{align}
for some norm $\norm{\cdot}$ on $\R^{n^d}$.

With regard to answering \cref{eq:basicq}, it is well-known that the following
dichotomy holds: while for ``most''\footnote{E.g., the set of singular matrices, as a subset of $\R^{n\times n}$, has Lebesgue measure zero.
By extension, any $\nu$-unfolding of
a tensor is a matrix in
$\R^{n^\nu\times n^{d-\nu}}$, and thus the
set of tensors with \emph{exact} low-rank representations
has Lebesgue measure zero.
Similar statements apply for
$\varepsilon$-approximations in
an appropriate sense, to be explained shortly.} tensors the rank $r_\varepsilon$ will
scale exponentially with $d$, in many practical applications $r_\varepsilon$
will scale at most polynomially with $d$.
For instance, in \cite{Hackbusch2005,Hackbusch2005a}, the authors show that if the tensor in question is
an analytic function of a sum of the independent variables
or it is a power of a discrete Laplacian operator, then, such tensors
can be very efficiently approximated -- without the curse of dimensionality.
In other words, the curse of dimensionality can be circumvented
provided one assumes \emph{enough} structure.

What is then \emph{enough} structure?
One can assume a highly structured expression for a tensor
$\bs\psi\in\R^{n^d}$ or a linear operator
$A:\R^{n^d}\rightarrow\R^{n^d}$, and
an approximation method tailored to that expression,
similarly as in \cite{Hackbusch2005,Hackbusch2005a,HACKBUSCH2007697}.
One can assume $\bs\psi$ is a solution to the
Poisson equation with low-rank data -- as in, e.g., \cite{DahmenTS}.
These assumptions are certainly sufficient.
If, however, one only assumes $\bs\psi$ belongs to a regularity
class -- such as a Sobolev space -- the curse of dimensionality
cannot be avoided\footnote{This is not surprising as Kolmogorov
showed that, for a fixed precision, Sobolev balls require an amount of information
exponential in the dimension $d$ to encode
\cite{Tikhomirov1993}, i.e., these sets are too large and
do not have enough structure to avoid the curse of dimensionality.}
(see, e.g., \cite{Schneider}).

The main \textbf{contribution} of this work is along those lines.
Our setting is that of a system described by a Hamiltonian $H$ -- a self-adjoint
densely defined differential operator;
and we consider the scaling of $r_\varepsilon$ of
the ground state projection\footnote{A.k.a. Ground state density matrix/operator.} (GSP)
-- the orthogonal projection onto the eigenspace of the smallest eigenvalue
of $H$.
This setting is frequently encountered in quantum mechanics.
Our assumptions involve $H$ and its spectrum,
i.e., in particular, we do not assume any specific expression for
the eigenfunctions of $H$.

Naturally, one can make stronger statements about approximability if one has
explicit representations of $H$,
and devise efficient computational methods tailored to $H$.
However, as we will argue in \Cref{rem:ass},
many of the assumptions made are sufficient
\emph{and} necessary\footnote{In the sense that, if one assumption
is not satisfied, there are examples of physical systems for which
the ranks will scale exponentially.}.
The assumptions have a physical interpretation and are thus
not merely technicalities of the mathematical formalism.
Moreover, we indicate shortly
(see \Cref{sec:thiswork}) that our approach offers the flexibility to choose
more general structures of $H$.
The framework of our results thus comes close
to walking the fine line of just
enough assumptions to guarantee a target
function or operator
\emph{does not} suffer from the curse of dimensionality.


\subsection{Tensor Networks and Hierarchical Ranks}
\begin{figure}[ht!]
    \begin{subfigure}{0.47\textwidth}
        \center
        \tikzset{every picture/.style={line width=0.75pt}} 

\begin{tikzpicture}[x=0.75pt,y=0.75pt,yscale=-1,xscale=1]

\draw    (279.25,117.75) -- (279.25,145.5) ;
\draw    (360.25,117.75) -- (360.25,145.5) ;
\draw    (279.25,117.75) -- (360.25,117.75) ;
\draw  [draw opacity=0][fill={rgb, 255:red, 74; green, 144; blue, 226 }  ,fill opacity=1 ] (275.5,117.75) .. controls (275.5,115.68) and (277.18,114) .. (279.25,114) .. controls (281.32,114) and (283,115.68) .. (283,117.75) .. controls (283,119.82) and (281.32,121.5) .. (279.25,121.5) .. controls (277.18,121.5) and (275.5,119.82) .. (275.5,117.75) -- cycle ;
\draw  [draw opacity=0][fill={rgb, 255:red, 74; green, 144; blue, 226 }  ,fill opacity=1 ] (356.5,117.75) .. controls (356.5,115.68) and (358.18,114) .. (360.25,114) .. controls (362.32,114) and (364,115.68) .. (364,117.75) .. controls (364,119.82) and (362.32,121.5) .. (360.25,121.5) .. controls (358.18,121.5) and (356.5,119.82) .. (356.5,117.75) -- cycle ;

\end{tikzpicture}
        \caption{TT corresponding to \cref{eq:tt1}
        with \emph{open} boundary conditions.}
    \end{subfigure}
    
    \begin{subfigure}{0.47\textwidth}
        \center
        \tikzset{every picture/.style={line width=0.75pt}} 

\begin{tikzpicture}[x=0.75pt,y=0.75pt,yscale=-1,xscale=1]

\draw    (441.25,117.75) -- (441.25,145.5) ;
\draw    (360.25,117.75) -- (441.25,117.75) ;
\draw    (279.25,117.75) -- (279.25,145.5) ;
\draw    (360.25,117.75) -- (360.25,145.5) ;
\draw    (279.25,117.75) -- (360.25,117.75) ;
\draw  [draw opacity=0][fill={rgb, 255:red, 74; green, 144; blue, 226 }  ,fill opacity=1 ] (275.5,117.75) .. controls (275.5,115.68) and (277.18,114) .. (279.25,114) .. controls (281.32,114) and (283,115.68) .. (283,117.75) .. controls (283,119.82) and (281.32,121.5) .. (279.25,121.5) .. controls (277.18,121.5) and (275.5,119.82) .. (275.5,117.75) -- cycle ;
\draw  [draw opacity=0][fill={rgb, 255:red, 74; green, 144; blue, 226 }  ,fill opacity=1 ] (356.5,117.75) .. controls (356.5,115.68) and (358.18,114) .. (360.25,114) .. controls (362.32,114) and (364,115.68) .. (364,117.75) .. controls (364,119.82) and (362.32,121.5) .. (360.25,121.5) .. controls (358.18,121.5) and (356.5,119.82) .. (356.5,117.75) -- cycle ;
\draw  [draw opacity=0][fill={rgb, 255:red, 74; green, 144; blue, 226 }  ,fill opacity=1 ] (437.5,117.75) .. controls (437.5,115.68) and (439.18,114) .. (441.25,114) .. controls (443.32,114) and (445,115.68) .. (445,117.75) .. controls (445,119.82) and (443.32,121.5) .. (441.25,121.5) .. controls (439.18,121.5) and (437.5,119.82) .. (437.5,117.75) -- cycle ;

\end{tikzpicture}
        \caption{TT corresponding to \cref{eq:tt2}
        with \emph{open} boundary conditions.}
    \end{subfigure}

    \begin{subfigure}{0.47\textwidth}
        \center
        \tikzset{every picture/.style={line width=0.75pt}} 

\begin{tikzpicture}[x=0.75pt,y=0.75pt,yscale=-1,xscale=1]

\draw    (219.25,117.75) .. controls (259.25,87.75) and (368.67,61.5) .. (458.5,117.75) ;
\draw    (381.25,117.75) -- (462.25,117.75) ;
\draw    (458.5,117.75) -- (458.5,145.5) ;
\draw    (381.25,117.75) -- (381.25,145.5) ;
\draw    (300.25,117.75) -- (381.25,117.75) ;
\draw    (219.25,117.75) -- (219.25,145.5) ;
\draw    (300.25,117.75) -- (300.25,145.5) ;
\draw    (219.25,117.75) -- (300.25,117.75) ;
\draw  [draw opacity=0][fill={rgb, 255:red, 74; green, 144; blue, 226 }  ,fill opacity=1 ] (215.5,117.75) .. controls (215.5,115.68) and (217.18,114) .. (219.25,114) .. controls (221.32,114) and (223,115.68) .. (223,117.75) .. controls (223,119.82) and (221.32,121.5) .. (219.25,121.5) .. controls (217.18,121.5) and (215.5,119.82) .. (215.5,117.75) -- cycle ;
\draw  [draw opacity=0][fill={rgb, 255:red, 74; green, 144; blue, 226 }  ,fill opacity=1 ] (296.5,117.75) .. controls (296.5,115.68) and (298.18,114) .. (300.25,114) .. controls (302.32,114) and (304,115.68) .. (304,117.75) .. controls (304,119.82) and (302.32,121.5) .. (300.25,121.5) .. controls (298.18,121.5) and (296.5,119.82) .. (296.5,117.75) -- cycle ;
\draw  [draw opacity=0][fill={rgb, 255:red, 74; green, 144; blue, 226 }  ,fill opacity=1 ] (377.5,117.75) .. controls (377.5,115.68) and (379.18,114) .. (381.25,114) .. controls (383.32,114) and (385,115.68) .. (385,117.75) .. controls (385,119.82) and (383.32,121.5) .. (381.25,121.5) .. controls (379.18,121.5) and (377.5,119.82) .. (377.5,117.75) -- cycle ;
\draw  [draw opacity=0][fill={rgb, 255:red, 74; green, 144; blue, 226 }  ,fill opacity=1 ] (454.75,117.75) .. controls (454.75,115.68) and (456.43,114) .. (458.5,114) .. controls (460.57,114) and (462.25,115.68) .. (462.25,117.75) .. controls (462.25,119.82) and (460.57,121.5) .. (458.5,121.5) .. controls (456.43,121.5) and (454.75,119.82) .. (454.75,117.75) -- cycle ;

\end{tikzpicture}
        \caption{TT corresponding to \cref{eq:tt3}
        with \emph{periodic} boundary conditions.}
    \end{subfigure}

    \begin{subfigure}{0.47\textwidth}
        \center
        \tikzset{every picture/.style={line width=0.75pt}} 

\begin{tikzpicture}[x=0.75pt,y=0.75pt,yscale=-1,xscale=1]

\draw    (225.25,185.75) -- (225.25,212.5) ;
\draw    (306.25,185.75) -- (306.25,212.5) ;
\draw    (387.25,185.75) -- (387.25,212.5) ;
\draw    (464.5,185.75) -- (464.5,212.5) ;
\draw    (425.85,131.75) -- (464.5,185.75) ;
\draw    (425.85,131.75) -- (387.25,185.75) ;
\draw    (265.75,131.75) -- (306.25,185.75) ;
\draw    (265.75,131.75) -- (225.25,185.75) ;
\draw    (346.25,77.75) -- (265.75,131.75) ;
\draw    (346.25,77.75) -- (425.85,131.75) ;
\draw  [draw opacity=0][fill={rgb, 255:red, 74; green, 144; blue, 226 }  ,fill opacity=1 ] (221.5,185.75) .. controls (221.5,183.68) and (223.18,182) .. (225.25,182) .. controls (227.32,182) and (229,183.68) .. (229,185.75) .. controls (229,187.82) and (227.32,189.5) .. (225.25,189.5) .. controls (223.18,189.5) and (221.5,187.82) .. (221.5,185.75) -- cycle ;
\draw  [draw opacity=0][fill={rgb, 255:red, 74; green, 144; blue, 226 }  ,fill opacity=1 ] (302.5,185.75) .. controls (302.5,183.68) and (304.18,182) .. (306.25,182) .. controls (308.32,182) and (310,183.68) .. (310,185.75) .. controls (310,187.82) and (308.32,189.5) .. (306.25,189.5) .. controls (304.18,189.5) and (302.5,187.82) .. (302.5,185.75) -- cycle ;
\draw  [draw opacity=0][fill={rgb, 255:red, 74; green, 144; blue, 226 }  ,fill opacity=1 ] (383.5,185.75) .. controls (383.5,183.68) and (385.18,182) .. (387.25,182) .. controls (389.32,182) and (391,183.68) .. (391,185.75) .. controls (391,187.82) and (389.32,189.5) .. (387.25,189.5) .. controls (385.18,189.5) and (383.5,187.82) .. (383.5,185.75) -- cycle ;
\draw  [draw opacity=0][fill={rgb, 255:red, 74; green, 144; blue, 226 }  ,fill opacity=1 ] (460.75,185.75) .. controls (460.75,183.68) and (462.43,182) .. (464.5,182) .. controls (466.57,182) and (468.25,183.68) .. (468.25,185.75) .. controls (468.25,187.82) and (466.57,189.5) .. (464.5,189.5) .. controls (462.43,189.5) and (460.75,187.82) .. (460.75,185.75) -- cycle ;
\draw  [draw opacity=0][fill={rgb, 255:red, 74; green, 144; blue, 226 }  ,fill opacity=1 ] (262,131.75) .. controls (262,129.68) and (263.68,128) .. (265.75,128) .. controls (267.82,128) and (269.5,129.68) .. (269.5,131.75) .. controls (269.5,133.82) and (267.82,135.5) .. (265.75,135.5) .. controls (263.68,135.5) and (262,133.82) .. (262,131.75) -- cycle ;
\draw  [draw opacity=0][fill={rgb, 255:red, 74; green, 144; blue, 226 }  ,fill opacity=1 ] (422.1,131.75) .. controls (422.1,129.68) and (423.78,128) .. (425.85,128) .. controls (427.92,128) and (429.6,129.68) .. (429.6,131.75) .. controls (429.6,133.82) and (427.92,135.5) .. (425.85,135.5) .. controls (423.78,135.5) and (422.1,133.82) .. (422.1,131.75) -- cycle ;
\draw  [draw opacity=0][fill={rgb, 255:red, 74; green, 144; blue, 226 }  ,fill opacity=1 ] (342.5,77.75) .. controls (342.5,75.68) and (344.18,74) .. (346.25,74) .. controls (348.32,74) and (350,75.68) .. (350,77.75) .. controls (350,79.82) and (348.32,81.5) .. (346.25,81.5) .. controls (344.18,81.5) and (342.5,79.82) .. (342.5,77.75) -- cycle ;

\end{tikzpicture}
        \caption{A tree-based TN.}\label{fig:treebased}
    \end{subfigure}
    
    \begin{subfigure}{0.47\textwidth}
        \center
        \tikzset{every picture/.style={line width=0.75pt}} 

\begin{tikzpicture}[x=0.75pt,y=0.75pt,yscale=-1,xscale=1]

\draw    (262.25,77.75) .. controls (302.25,47.75) and (411.67,73.5) .. (427.67,107.11) ;
\draw    (262.25,51) -- (262.25,77.75) ;
\draw    (460,48.03) -- (460,74.78) ;
\draw    (492.33,107.11) -- (520,107.11) ;
\draw    (460,139.45) -- (460,166.2) ;
\draw    (403.03,131.75) -- (427.67,107.11) ;
\draw    (376.75,185.75) -- (376.75,212.5) ;
\draw    (315.43,185.75) -- (315.43,212.5) ;
\draw    (141.25,185.75) -- (222.25,185.75) ;
\draw   (341.71,131.75) -- (403.03,131.75) -- (376.75,185.75) -- (315.43,185.75) -- cycle ;
\draw    (141.25,185.75) -- (141.25,212.5) ;
\draw    (222.25,185.75) -- (222.25,212.5) ;
\draw    (181.75,131.75) -- (222.25,185.75) ;
\draw    (181.75,131.75) -- (141.25,185.75) ;
\draw    (262.25,77.75) -- (181.75,131.75) ;
\draw    (262.25,77.75) -- (341.85,131.75) ;
\draw  [draw opacity=0][fill={rgb, 255:red, 74; green, 144; blue, 226 }  ,fill opacity=1 ] (137.5,185.75) .. controls (137.5,183.68) and (139.18,182) .. (141.25,182) .. controls (143.32,182) and (145,183.68) .. (145,185.75) .. controls (145,187.82) and (143.32,189.5) .. (141.25,189.5) .. controls (139.18,189.5) and (137.5,187.82) .. (137.5,185.75) -- cycle ;
\draw  [draw opacity=0][fill={rgb, 255:red, 74; green, 144; blue, 226 }  ,fill opacity=1 ] (218.5,185.75) .. controls (218.5,183.68) and (220.18,182) .. (222.25,182) .. controls (224.32,182) and (226,183.68) .. (226,185.75) .. controls (226,187.82) and (224.32,189.5) .. (222.25,189.5) .. controls (220.18,189.5) and (218.5,187.82) .. (218.5,185.75) -- cycle ;
\draw  [draw opacity=0][fill={rgb, 255:red, 74; green, 144; blue, 226 }  ,fill opacity=1 ] (311.68,185.75) .. controls (311.68,183.68) and (313.36,182) .. (315.43,182) .. controls (317.5,182) and (319.18,183.68) .. (319.18,185.75) .. controls (319.18,187.82) and (317.5,189.5) .. (315.43,189.5) .. controls (313.36,189.5) and (311.68,187.82) .. (311.68,185.75) -- cycle ;
\draw  [draw opacity=0][fill={rgb, 255:red, 74; green, 144; blue, 226 }  ,fill opacity=1 ] (373,185.75) .. controls (373,183.68) and (374.68,182) .. (376.75,182) .. controls (378.82,182) and (380.5,183.68) .. (380.5,185.75) .. controls (380.5,187.82) and (378.82,189.5) .. (376.75,189.5) .. controls (374.68,189.5) and (373,187.82) .. (373,185.75) -- cycle ;
\draw  [draw opacity=0][fill={rgb, 255:red, 74; green, 144; blue, 226 }  ,fill opacity=1 ] (178,131.75) .. controls (178,129.68) and (179.68,128) .. (181.75,128) .. controls (183.82,128) and (185.5,129.68) .. (185.5,131.75) .. controls (185.5,133.82) and (183.82,135.5) .. (181.75,135.5) .. controls (179.68,135.5) and (178,133.82) .. (178,131.75) -- cycle ;
\draw  [draw opacity=0][fill={rgb, 255:red, 74; green, 144; blue, 226 }  ,fill opacity=1 ] (338.1,131.75) .. controls (338.1,129.68) and (339.78,128) .. (341.85,128) .. controls (343.92,128) and (345.6,129.68) .. (345.6,131.75) .. controls (345.6,133.82) and (343.92,135.5) .. (341.85,135.5) .. controls (339.78,135.5) and (338.1,133.82) .. (338.1,131.75) -- cycle ;
\draw  [draw opacity=0][fill={rgb, 255:red, 74; green, 144; blue, 226 }  ,fill opacity=1 ] (258.5,77.75) .. controls (258.5,75.68) and (260.18,74) .. (262.25,74) .. controls (264.32,74) and (266,75.68) .. (266,77.75) .. controls (266,79.82) and (264.32,81.5) .. (262.25,81.5) .. controls (260.18,81.5) and (258.5,79.82) .. (258.5,77.75) -- cycle ;
\draw  [draw opacity=0][fill={rgb, 255:red, 74; green, 144; blue, 226 }  ,fill opacity=1 ] (399.28,131.75) .. controls (399.28,129.68) and (400.96,128) .. (403.03,128) .. controls (405.1,128) and (406.78,129.68) .. (406.78,131.75) .. controls (406.78,133.82) and (405.1,135.5) .. (403.03,135.5) .. controls (400.96,135.5) and (399.28,133.82) .. (399.28,131.75) -- cycle ;
\draw   (427.67,107.11) .. controls (427.67,89.26) and (442.14,74.78) .. (460,74.78) .. controls (477.86,74.78) and (492.33,89.26) .. (492.33,107.11) .. controls (492.33,124.97) and (477.86,139.45) .. (460,139.45) .. controls (442.14,139.45) and (427.67,124.97) .. (427.67,107.11) -- cycle ;
\draw  [draw opacity=0][fill={rgb, 255:red, 74; green, 144; blue, 226 }  ,fill opacity=1 ] (423.92,107.11) .. controls (423.92,105.04) and (425.6,103.36) .. (427.67,103.36) .. controls (429.74,103.36) and (431.42,105.04) .. (431.42,107.11) .. controls (431.42,109.18) and (429.74,110.86) .. (427.67,110.86) .. controls (425.6,110.86) and (423.92,109.18) .. (423.92,107.11) -- cycle ;
\draw  [draw opacity=0][fill={rgb, 255:red, 74; green, 144; blue, 226 }  ,fill opacity=1 ] (456.25,139.45) .. controls (456.25,137.38) and (457.93,135.7) .. (460,135.7) .. controls (462.07,135.7) and (463.75,137.38) .. (463.75,139.45) .. controls (463.75,141.52) and (462.07,143.2) .. (460,143.2) .. controls (457.93,143.2) and (456.25,141.52) .. (456.25,139.45) -- cycle ;
\draw  [draw opacity=0][fill={rgb, 255:red, 74; green, 144; blue, 226 }  ,fill opacity=1 ] (488.58,107.11) .. controls (488.58,105.04) and (490.26,103.36) .. (492.33,103.36) .. controls (494.4,103.36) and (496.08,105.04) .. (496.08,107.11) .. controls (496.08,109.18) and (494.4,110.86) .. (492.33,110.86) .. controls (490.26,110.86) and (488.58,109.18) .. (488.58,107.11) -- cycle ;
\draw  [draw opacity=0][fill={rgb, 255:red, 74; green, 144; blue, 226 }  ,fill opacity=1 ] (456.25,74.78) .. controls (456.25,72.71) and (457.93,71.03) .. (460,71.03) .. controls (462.07,71.03) and (463.75,72.71) .. (463.75,74.78) .. controls (463.75,76.85) and (462.07,78.53) .. (460,78.53) .. controls (457.93,78.53) and (456.25,76.85) .. (456.25,74.78) -- cycle ;

\end{tikzpicture}
        \caption{A general TN.}
    \end{subfigure}

    \caption{Examples of tensor networks. A vertex corresponds to
    a \emph{tensor core} $U^j$, connected edges correspond to indices
    in $U^j$ that are summed over or \emph{contracted},
    free edges correspond to free input indices/variables.}
    \label{fig:tns}
\end{figure}
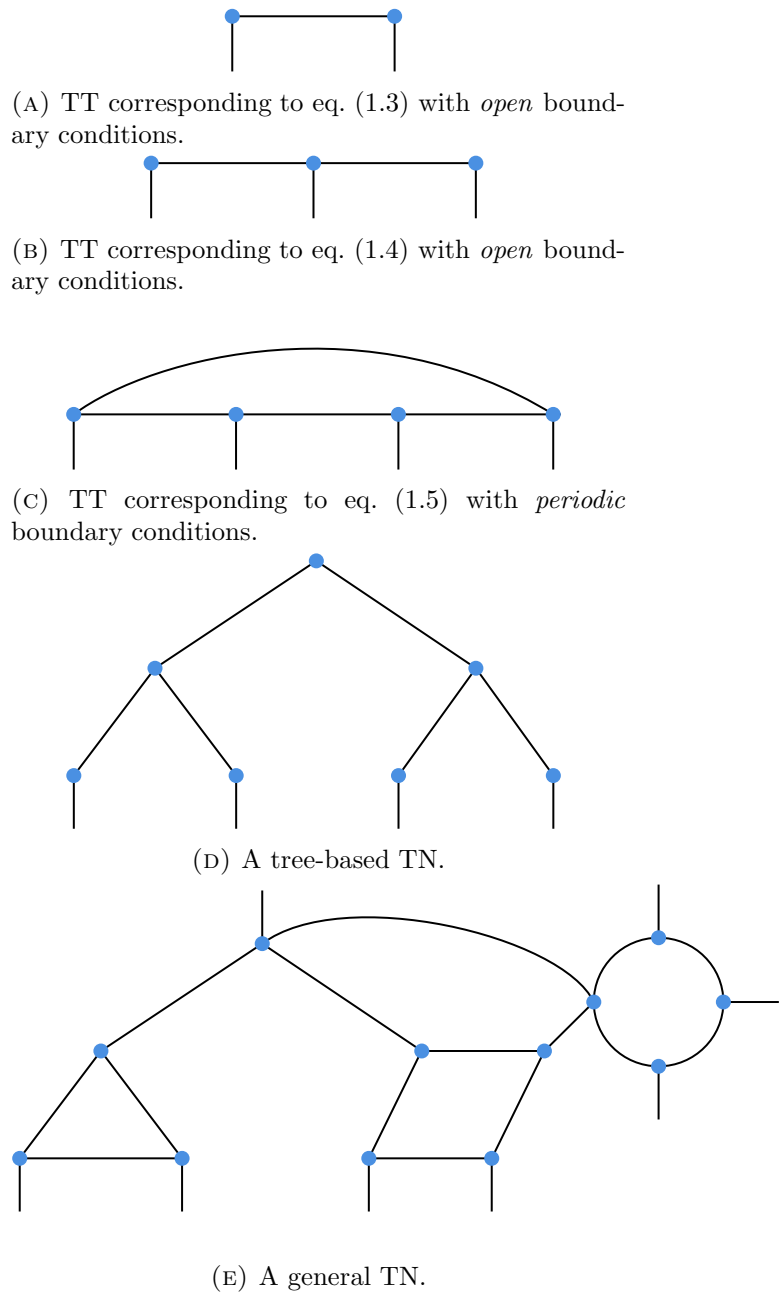

At this point we have to specify that our results do not concern
approximation as in \cref{eq:cp} but rather approximation
with Tensor Networks (TNs).
Approximations as in \cref{eq:cp} are frequently
referred to as $r$-term approximation or
Canonical Polyadic (CP) decomposition.
Though CP certainly has its uses -- most notably in chemometrics
\cite{HARSHMAN199439} --
this format is not flexible\footnote{Apart from also being numerically
unstable \cite{RePEc}.} enough for efficiently
approximating more complex tensor product structures.
The type of format we need for this work is a Tensor Train (TT) \cite{Oseledets2009},
a.k.a.\ a Matrix Product State (MPS) \cite{Kl_mper_1993} -- a particular
type of a TN. The TT of an array $\bs\psi\in\R^{n\times\cdots\times n}$
has the form (for $d=2,3,4$)
\begin{align}
    \bs\psi(i_1,i_2)&=\sum_{k_1=1}^{r_1}U^1(i_1,k_1)U^2(k_1,i_2),\label{eq:tt1}\\
    \bs\psi(i_1,i_2,i_3)&=\sum_{k_1=1}^{r_1}
    \sum_{k_2=1}^{r_2}U^1(i_1,k_1)U^2(k_1,i_2,k_2)U^3(k_2,i_3),\label{eq:tt2}\\
    \bs\psi(i_1,i_2,i_3,i_4)&=\sum_{k_0=1}^{r_0}\sum_{k_1=1}^{r_1}
    \sum_{k_2=1}^{r_2}
    \sum_{k_3=1}^{r_3}U^1(k_0,i_1,k_1)U^2(k_1,i_2,k_2)U^3(k_2,i_3,k_3)
    U^4(k_3,i_4,k_0),\label{eq:tt3}
\end{align}
with order-2 or order-3 tensors $U^j\in\R^{r_{j-1}\times n\times r_j}$.
It is convenient to visualize and manipulate
tensor networks via tensor diagrams\footnote{A.k.a.\
Penrose graphical notation.}, see \Cref{fig:tns}.
As we will indicate in  \Cref{sec:thiswork},
our work can be in principle extended to more general networks.

Our main result concerns \emph{hierarchical} or \emph{multi-linear} ranks.
These are ranks of bi-partite cuts,
i.e., for any $\alpha\subset\{1,\ldots,d\}$ and
$\alpha^c:=\{1,\ldots,d\}\setminus\alpha$,
the hierarchical rank $r_\alpha=r_{\alpha^c}$ is the smallest
$r_\alpha\in\N_0$ such that
\begin{align}\label{eq:hranks}
    \bs\psi=\sum_{k_\alpha=1}^{r_\alpha}\bs v_{k_\alpha}^{\alpha}\otimes
    \bs w_{k_\alpha}^{\alpha^c},\quad
    \bs v^\alpha_{k_\alpha}\in \R^{\otimes\#\alpha},\;
    \bs w^{\alpha^c}_{k_\alpha}\in \R^{\otimes\#\alpha^c}.
\end{align}
These ranks essentially determine the representation and computational
complexity of a TN.
The case of TTs corresponds to bi-partite cuts of the form $\alpha=\{1,
\ldots,j\}$ for $j=1,\ldots,d-1$. In this case we abbreviate $r_j:=
r_{\{1,\ldots,j\}}$.
From \Cref{fig:tns} it is not difficult to see that the representation complexity of
tree tensor networks
scales roughly linearly (or log-linearly) with $d$ and
polynomially in ranks $\mc O(r_\alpha^p)$, where
$p$ is the number of edges connecting a vertex.

Notice the one-way relationship between the hierarchical ranks
$r_\alpha$ and the tensor rank $r$ from \cref{eq:cp}:
$r_\alpha\leq r$ for any $\alpha$, while the reverse inequality is not true
in general.
Loosely speaking, \cref{eq:cp} demands simultaneous separability
w.r.t.\ \emph{all} coordinate directions,
while \cref{eq:hranks} demands separability only w.r.t.\
bi-partite cuts. And thus one can have TNs with small ranks
$r_\alpha$ but exponentially growing tensor rank $r=\mc O(2^d)$.
On the other hand, once again assuming, e.g., only
classical Sobolev regularity,
$r_\alpha$ will scale exponentially with $d$ as well \cite{Schneider}.


\subsection{Contribution}
To put this work into context, we summarize what is known about
low-rank approximation with TNs and how our work
contributes to that body of knowledge.
We do not intend to give a complete overview of all the research on
TN-approximation -- a lot has been done in the past decades.
Instead, we highlight a few works that will serve as representatives
of a general flavor of results
relevant to us.
We apologize in advance if some contributors
feel left out by this summary.

\subsubsection*{Previous Work}
There are various results that construct efficient approximations to functions and operators
with explicit or highly structured expressions,
mostly in the context of numerical approximation
and finite-dimensional Hilbert spaces as in, e.g., \cite{Hackbusch2005,
Hackbusch2005a,QTTCayley}.
Techniques such as exponential sum approximation or
sinc quadrature are key for these methods.

Some results bound the ranks of an iterative method on
finite-dimensional Hilbert spaces, see, e.g., \cite{HBetal2012,Khoromskij2009TensorStructuredPA}.
There are a few results for infinite-dimensional Hilbert spaces as in
\cite{DahmenTS,BD,Ali2020,Bachmayr2014,BDL2}.
In \cite{BD,Ali2020,Bachmayr2014,BDL2}, the authors consider adaptive numerical
methods for PDEs and show that,
if the \emph{exact} solution is low-rank approximable,
then the same holds for the numerical approximation.
In \cite{DahmenTS},
the authors show for the Poisson equation that the
exact solution has
similar low-rank approximability as the right-hand-side data.

In \cite{Schneider}, the authors compute the optimal rank scaling for periodic functions in
Sobolev and mixed Sobolev spaces,
with the curse of dimensionality present in both cases, though
much milder for mixed Sobolev.
In \cite{Harbrecht}, the authors used weighted Sobolev spaces where, loosely speaking,
the weighting ensures that for
$d\rightarrow\infty$ the additional
dimensions have smaller and smaller impact on the
derivative of the function.
For these spaces the optimal
hierarchical ranks scale without the curse of dimensionality.

Note that classical regularity is certainly important
for a complete approximation scheme,
as even if one has a rank-one function/operator, in general,
one still has to approximate the low-dimensional components.
If the latter is done via, e.g., piece-wise
polynomials, then its efficiency depends
on the degree of classical Sobolev/Besov/analytic regularity
of said low-dimensional components.

However, if one is interested \emph{solely}
in the question of rank scaling -- as we are in this work
-- then this property should be independent
of any notion of smoothness.
An analytic function of the sum of
$d$ variables can be certainly approximated with ranks
scaling without
the curse of dimensionality -- for any TN.
But a tensor product of delta distributions has rank one in any representation --
despite the absence of any notion of classical regularity.
Hence, we emphasize that the techniques used in this work do not rely
on any classical notion of regularity.

\subsubsection{This Work}\label{sec:thiswork}\leavevmode
\begin{enumerate}[label=(\roman*)]
    \item We provide a rigorous, discretization-
    and (classical) regularity-independent proof of how
    a TN that is adapted to the Hamiltonian structure avoids the curse
    of dimensionality.
    
    \item Many of the required assumptions are necessary, have a physical
    interpretation and establish a direct link
    to the presence or absence of the
    curse of dimensionality.
    They provide further insight into when and why
    said curse appears, which goes beyond
    the mere intuition
    ``Hamiltonian structure $\approx$ TN structure''.
    
    \item The utilized technique is at least as valuable as the result itself.
    It demonstrates how similar results can be derived for more
    complex TNs, evolution equations and PDEs. E.g.,
    the approach can be extended to Hamiltonians
    with $k$-neighbor or long-range interactions that decay
    sufficiently fast, see also \cite{RST}.
    It can be extended to more general networks
    using Lieb-Robinson bounds from, e.g.,
    \cite{LBInfDim}.
    It provides an alternative to exponential sum approximations
    in the case of a unitary evolution via
    $\exp{-\i Ht}$, without the usual decay of eigenvalues.
    Finally, an additional benefit of this technique
    is that it provides valuable physical insight.
    
    \item Our goal is to identify the right framework for describing low-rank
    approximability or, more generally,
    approximation without the curse of dimensionality.
    We assume just enough to guarantee approximability,
    while still leaving flexibility in the setting where
    possible. This contribution is an important step towards that goal.
\end{enumerate}


\subsection*{Outline}
In \Cref{sec:result},
we introduce the concept of operator support,
our main result and discuss implications.
\Cref{sec:alaw} is dedicated to the proof
of the main result in \Cref{thm:obolor}.
The technical proofs are postponed to \Cref{app:proofs}.
\section{Main Result and Discussion}\label{sec:result}

Our main statement requires the concept of operator support.
Though our result concerns general (possibly infinite-dimensional)
Hilbert spaces,
for the purpose of introduction we mainly stick with finite
vectors and matrices.

Let $A\in\mc L\left((\R^n)^{\otimes d}\right)$ be a bounded linear operator
from $(\R^n)^{\otimes d}$ onto itself (an $n^d\times n^d$ matrix).
If for some $\alpha\subset\{1,\ldots,d\}$ and some $A_\alpha
\in\mc L\left((\R^n)^{\otimes\#\alpha}\right)$ we can write $A$ as
\begin{align*}
    A=A_\alpha\otimes\id_{\alpha^c},
\end{align*}
where $\id_{\alpha^c}:(\R^n)^{\otimes\#\alpha^c}\rightarrow
(\R^n)^{\otimes\#\alpha^c}$ is the identity, then we say
$A$ is \emph{supported} on the Hilbert space
$\H_\alpha:=(\R^n)^{\otimes\#\alpha}$,
or simply \emph{supported on $\alpha$},
or \emph{acts} on $\H_\alpha$.
We write in this case $\supp{A}:=\H_\alpha$
to denote the support of $A$.

For operators $A$, $B$, we write
$\dist(A,B)$ to denote the distance
between supports. E.g., suppose
$A,B\in\mc L(\mc H_1\otimes\ldots\otimes\mc H_d)$,
$\supp{A}= \mc H_2\otimes \mc H_3$ and
$\supp{B}=\mc H_3\otimes\mc H_4$.
Then, $\dist(A,B)=0$. If
$\supp{A}= \mc H_1\otimes \mc H_2$ and
$\supp{B}=\mc H_4\otimes\mc H_5$, then
$\dist(A,B)=2$, and so on.

The relevance of operator support is as follows:
if $A$ is supported on $\alpha$, then,
since the identity $\id_{\alpha^c}$ is a tensor product rank-one
operator,
we know that a low-rank approximation of
$A_\alpha$ w.r.t.\ any cut $\gamma\subset\alpha$
grows at most exponentially in $\#\alpha/2$ -- and not $d/2$.
Thus, if $\#\alpha$ can be guaranteed to be small
independently of $d$, we have decoupled
the global (large) dimension $d$ from the
effective dimension $\#\alpha$ of $A$, i.e.,
the dimension of the space on which $A$ acts non-trivially.

The core technique of our proof is the quantification
of support growth due to Lieb-Robinson \cite{LBOriginal}:
one starts with a rank-one operator, which is then
smoothly evolved/continued into the target state.
In the process of said continuation\footnote{Which is often referred to
as \emph{adiabatic continuation} in physics.},
the support of the operator grows at a speed determined
by the properties of the system Hamiltonian $H$.

The main ideas behind this approach are based on the
physical principle of holography: the informational content/entropy of
many physical systems is proportional to its area rather than volume.
Specifically, we build on the quantum entanglement area laws research
of \cite{Hastings,QA} and Lieb-Robinson bounds from
\cite{LRComm,LBInfDim,NearComm}.
A more detailed exposition of our work including connections to
entanglement entropy area laws can be found in \cite{ali19}.

\subsection{Main Result}
Now that we have introduced all the necessary concepts for this work,
we state our main result. We consider separable complex-valued Hilbert spaces
$\mc H_j$, $j=1,\ldots, d$, and their tensor product space
$\H:=\bigotimes_{j=1}^d\mc H_j$. Here, we consider general topological
tensor spaces $\H$ equipped and completed with the
canonical Hilbert space inner product.
If at least $d-1$ of the Hilbert spaces are finite-dimensional, then
the distinction between algebraic and topological tensor product
spaces is irrelevant. We also introduce the short-hand notation
\begin{align*}
    \bs{\mc H}_{i,j}:=\bigotimes_{\nu=i}^j\mc H_\nu,
\end{align*}
for the partial tensor product, once again equipped with
the canonical inner product.

Our Hamiltonian is a densely defined self-adjoint operator
$H:\dom{\H}\rightarrow\H$.
Our first primary assumption is on the locality
structure of $H$
\begin{align}\label{eq:introham}
    H=\sum_{j=1}^{d-1} H_{j,j+1},
\end{align}
where each $H_{j,j+1}:\dom{\H}\rightarrow\H$ is a Hamiltonian
in its own right, supported on $\bs{\mc H}_{j,j+1}$.

This is (loosely) to be interpreted as follows: for a multi-partite (e.g., many-body) system
with $d$ constituents,
$H_{j,j+1}$ models the joint energy of subsystems $j$ and $j+1$,
while the total energy of the entire system is described by $H$
and it is the sum of the subsystem energies of pairs $(j,j+1)$.
Such a description of a physical system is a reasonable approximation
if the many-body system has only pairwise interactions -- as is the
case for some systems of condensed matter,
see \cite{Carr2011}.

The particular local structure of $H$ in \cref{eq:introham}
fits the TT format. However, since Lieb-Robinson bounds hold for more general
graphs (see, e.g., \cite{LBInfDim}), one could extend our results to more general
TNs, using similar ideas but a more cumbersome proof.

The format of writing $H$ as in \cref{eq:introham} might seem
unorthodox if one is more accustomed to
\begin{align*}
    H:=K+V=-\frac{1}{2}\Delta+ V,
\end{align*}
where $\Delta$ is the Laplacian modeling the kinetic part and $V$
is some potential. Note, however, that if we can
write the potential as a superposition of local potential terms
$V=\sum_{j=1}^{d-1}V_{j,j+1}$,
where each local potential $V_{j,j+1}$
acts only on $\bs{\mc H}_{j,j+1}$, then
this provides a representation as in \cref{eq:introham},
since $-\Delta=-\sum_{j=1}^{d}\Delta_j$,
where $\Delta_j$ is a Laplacian acting on $\mc H_j$.

Our second primary assumption concerns the spectrum of $H$:
we suppose $H$ has a smallest eigenvalue $\lambda_0$,
which we set w.l.o.g.\ to $\lambda_0=0$, and that the spectral (energy)
gap $\Delta E:=\lambda_1-\lambda_0>0$ is non-vanishing for
$d\rightarrow\infty$, i.e., $\Delta E\not\rightarrow 0$
for $d\rightarrow\infty$.
This assumption separates critical systems with topological
degeneracy ($\Delta E\rightarrow 0$) from
non-critical systems without degeneracy ($\Delta E\not\rightarrow 0$):
for the former it is known, in general, that $r_j$ will scale exponentially with $d$.
Quantum criticality and the related exponential rank growth are a physical
property of the underlying system described by $H$ that
are in fact necessary to explain certain phenomena,
such as, e.g., superconductivity (see \cite{Carr2011}).

Finally, let the eigenspace
of $H$ corresponding to $\lambda_0$ be finite-dimensional
and $\{\bs\psi_0^k\}_{k=1}^N$ be an orthonormal basis of
said eigenspace. Then, we call the
(normalized\footnote{In principle,
one can consider any weighted average of
the one-dimensional projections
$ \inp{\cdot}{\bs\psi_0^k}_{\H}\bs\psi_0^k$.} trace-class) projection
\begin{align*}
    P_0=\frac{1}{N}\sum_{k=1}^N \inp{\cdot}{\bs\psi_0^k}_{\H}\bs\psi_0^k
    \in\mc L(\H),
\end{align*}
a \emph{ground state density operator} or
a \emph{ground state projection (GSP)}.

Then, we have
\begin{result}[\Cref{thm:obolor}]
    For any $j=1,\ldots,d-1$ and any $l=1,\ldots,\min\{j-1,d-j\}$, there exist
    bounded, positive, self-adjoint operators
    $B=B(j,l)$, $L=L(j,l)$ and $R=R(j,l)$ such that
    \begin{enumerate}[label=(\roman*)]
        \item $L$ and $R$ are supported on
        $\bs{\mc H}_{1,j}$ and $\bs{\mc H}_{j+1,d}$, respectively;
        \item $B$ is supported on $\bs{\mc H}_{j-l,j+l}$;
        \item the operator norms are bounded by unity
        $\norm{B}[\mc L]$, $\norm{L}[\mc L]$,
        $\norm{R}[\mc L]\leq 1$;
    \end{enumerate}
    such that we have the error estimate
    \begin{align}\label{eq:intromain}
        \norm{P_0-BLR}[\mc L]\leq C\exp{-cl},
    \end{align}
    where the constants $C,c>0$ depend only
    on $\Delta E$ and not $d$, $j$ or $l$.
\end{result}


\subsection{Implications}
The result of the above theorem provides a rigorous
and discretization-independent justification
of the well-accepted intuition that TNs that
\emph{fit} the (interaction) structure
of the problem -- in this case Hamiltonian
-- do not incur the curse of dimensionality.

\subsubsection{Hierarchical Ranks}
If for any $j=1,\ldots,d-1$, $P_0\approx B(j,l)L(j,l)R(j,l)$,
then the scaling of TT-ranks $r_j$ corresponding to bi-partite cuts
$\H:=\bs{\mc H}_{1,j}\otimes\bs{\mc H}_{j+1,d}$ can be
inferred from $BLR$. Note that the product $LR$ has rank-one
w.r.t.\ this bi-partite cut. Thus, a rank-$r_j$ approximation to $B$
is also a rank-$r_j$ approximation to $P_0$.

If the Hilbert spaces $\mc H_j$ are finite-dimensional with $\dim(\mc H_j)=n$,
then $B$ can be represented \emph{exactly}
with at most $r_j\leq n^l$, i.e., independently of $d$.
More generally, one would have to approximate $B$ as well.
Note that -- whatever the scaling of the ranks of an approximation
to $B$ is -- it will always be independent of $d$
\emph{as long as} the constants $C,c>0$ are independent of $d$.
We elaborate with an example.

Suppose $\mc H_j=L^2([0,1]^3)$ or $\mc H_j=L^2(\R^3)$.
As we will see in the proof of \Cref{thm:obolor},
$B$ is constructed from parts of the spectrum of $H$ and,
in many cases (see, e.g., \cite{Harry}), the ground state has some
Sobolev regularity $s>0$.
Recall that
\begin{align*}
    B=\hat{B}_l\otimes\id,\quad
    \hat{B}_l\in\mc L(\bs{\mc H}_{j-l,j+l}),\;
    \id\in\mc L(\bs{\mc H}_{1,j-l-1}\otimes\bs{\mc H}_{j+l+1,d}),
\end{align*}
and given Sobolev regularity of the kernel of $\hat{B}_l$,
\begin{align*}
    \norm{\hat{B}_l-\sum_{k=1}^{r_j}U_k\otimes V_k}[\mc L]
    =\mc O(r_j^{-s/l}),
    \quad U_k\in\mc L(\bs{\mc H}_{j-l,j}),\;
    V_k\in\mc L(\bs{\mc H}_{j+1,j+l}).
\end{align*}
Thus, for a fixed approximation accuracy $\varepsilon>0$,
from \cref{eq:intromain}
we first choose
$l:=\left\lceil1/c\ln(C/\varepsilon)\right\rceil$,
and hence
\begin{align*}
    r_j^{-s/l}\leq\varepsilon\quad
    \Leftrightarrow\quad
    r_j\geq \varepsilon^{-l/s}=
    \varepsilon^{(1/sc)\ln(\varepsilon/C)}.
\end{align*}
That is, exponential
dependence on $d$ can only be re-introduced
through $c$ or $C$.
E.g., if $c\sim 1/d$ or $C\sim 2^d$.
This in turn depends on $\Delta E$, i.e.,
the curse of dimensionality is now determined by the
spectral gap $\Delta E$ and, if the latter is bounded
from below independently of $d$,
then the curse of dimensionality can be (at least asymptotically) avoided.

The bound \cref{eq:intromain} decouples the dimension $d$ from
the scaling of hierarchical ranks or, alternatively,
shifts this dependence onto $C,c$.
The physical relationship between ranks $r_j$ and
constants $C,c$ is that both reflect the entanglement
properties of the ground state,
which in turn depend on the spectral gap,
local interaction length (in our case equal to two)
and interaction strength (see \Cref{ass:op}).
For gapped Hamiltonians, with fixed interaction length and strength,
the entanglement of the different subsystems saturates for $d\rightarrow\infty$:
there exist finite $C,c>0$, independent of $d$.
Consequently, ranks saturate as well.

\subsubsection{Evolution Equations}
A similar technique as in the proof of
\cref{eq:intromain} can be used to analyze evolution equations,
such as the Schrödinger equation (in Planck units)
\begin{align*}
    \i\frac{\d }{\d t}\bs\psi(t)=H\bs\psi(t),\quad
    \bs\psi(0)=\bs\psi_0.
\end{align*}
The density operator $P(t):=
\inp{\cdot}{\bs\psi(t)}_{\H}\bs\psi(t)$ evolves as
\begin{align*}
    P(t)=\exp{\i Ht}P(0)\exp{-\i Ht}.
\end{align*}
Starting with, e.g., a product state
$\bs\psi_0=\psi_0^1\otimes\ldots\otimes\psi_0^d$, similarly to \cref{eq:intromain}, one could show for any
$j=1,\ldots,d-1$
\begin{align*}
    P(t)\approx B(t,j,l)L(t,j,l)R(t,j,l).
\end{align*}
The support of $B(t,j,l)$ grows with $t\rightarrow\infty$ such that,
for small times $t>0$, $P(t)$ has small hierarchical ranks
independent of $d$.

Note the similarity to approximation techniques as in \cite{Hackbusch2005,Hackbusch2005a} or \cite{QTTCayley},
exploiting matrix exponentials with subsequent exponential sum approximations.
Unlike in the case of exponentials as
$\exp{-At}$ with a positive definite $A$ -- as for, e.g., parabolic evolution equations
-- the \emph{unitary} evolution via $\exp{-\i Ht}$ has no damping effect, i.e.,
the eigenvalues do not decay for $t\rightarrow\infty$.

For the case of unitary evolution,
in \cite{QTTCayley} the authors employ a quantized TT approximation based on
a series expansion of $\exp{-\i Ht}$ via they Cayley transform.
This guarantees highly efficient approximations
with low QTT-ranks of $\bs\psi(t)$
for the \emph{semi-discrete} problem (discrete in space),
\emph{provided} the initial value is itself of low rank,
analytic and the application of the Cayley transform remains of low rank.
In contrast,
the technique for proving
\cref{eq:intromain} using Lieb-Robinson bounds relies
solely on operator support. Hence, the
bound does not require finite-dimensional Hilbert spaces
or any notion of classical regularity,
but solely those assumptions that have direct influence
on the TN-ranks:
\begin{enumerate}[label=(\roman*)]
    \item a TN that fits the Hamiltonian structure,
    \item a Hamiltonian that does not allow interactions to spread
    arbitrarily fast (local/short-range and gapped),
    \item an initial state with small TN-ranks,
\end{enumerate}

\subsubsection{General Equations}
It is a common technique in functional calculus to express
functions of an operator as an integral, e.g.,
\begin{align*}
    A^{-\sigma}&=\frac{1}{\Gamma(\sigma)}\int_0^t
    t^{\sigma-1}\exp{-tA}\d t,\quad\sigma>0,\\
    F(A)&=\frac{1}{2\pi\i}\int_{\Gamma}F(z)(z-A)^{-1}\d z,
\end{align*}
see again \cite{Hackbusch2005,Hackbusch2005a,QTTCayley}.
This can be used, e.g., to solve
or pre-condition partial differential equations (PDEs).
Since the proof of \cref{eq:intromain} relies
on integral representations of
$P_0$ as well,
combined with subsequent Lieb-Robinson bounds, we believe there is
potential to extend such results to PDEs.

\subsubsection{Pure Ground States}
In many cases of interest the GS of $H$ is unique or \emph{pure},
i.e., $P_0=\inp{\cdot}{\bs\psi_0}_{\H}\bs\psi_0$.
In this case \cref{eq:intromain} does not immediately guarantee
low-rank approximability of $\bs\psi_0$.
To infer approximability of $\bs\psi_0$
one would have to show \cref{eq:intromain}
with the operator norm replaced by the trace norm.
Or, alternatively, one could extend the technique
of \cite{PhysRevB.85.195145,AradSubExp} to infinite-dimensional Hilbert spaces,
though this is not trivial.
In \cite{PhysRevB.85.195145,AradSubExp}, the authors use an approximate GSP --
similar to $BLR$ from \cref{eq:intromain}
but with stronger properties -- to conclude
that the pure GS $\bs\psi_0$ has non-exponentially growing TT/MPS-ranks. 
\section{Proof of \eqref{eq:intromain}}\label{sec:alaw}
We assume the Hamiltonian operator satisfies the following properties.

\begin{assumption}\label{ass:op}
    Let $H:\dom{H}\rightarrow\X$ be a densely defined
    self-adjoint (possibly unbounded) operator.
    \begin{enumerate}
        \item\label{ass:local} (Locality).
        We assume $H$ can be decomposed as
        \begin{align*}
            H=\sum_{j=1}^{d-1}H_{j,j+1},
        \end{align*}
        where each $H_{j,j+1}$ is supported on
        $\X_{j,j+1}$.

    \item\label{ass:gap} (Gap).
        We assume the spectrum is bounded from below
        with a non-vanishing spectral gap
        \begin{align}\label{eq:spgap}
            \Delta E:=\lambda_1-\lambda_0>0,
        \end{align}
        the GS eigenspace is finite-dimensional with
        dimension $N$ and
        the GSP is defined as
        \begin{align*}
            P_0 :=\frac{1}{N}\sum_{k=1}^\infty\inp{\cdot}{\bs\psi_0^k}_{\H}
            \bs\psi_0^k.
        \end{align*}

    \item\label{ass:int} (Finite Interaction Strength).
        We assume for all $1\leq j\leq d-1$, $H_{j,j+1}=H_j+H_{j+1}+
        \Phi_{j,j+1}$, where $H_j$ and $H_{j+1}$ are possibly unbounded
        operators supported on $\X_j$ and $\X_{j+1}$, respectively,
        and $\Phi_{j,j+1}$ is a densely defined uniformly bounded
        operator\footnote{That models interactions between particles
        $j$ and $j+1$.}
        supported on $\X_{j,j+1}$. I.e., there exists a constant $J$
        such that
        \begin{align}\label{eq:intstr}
            \norm{\Phi_{j,j+1}}[\mc L]\leq J,
        \end{align}
        for all $1\leq j\leq d-1$.

    \item\label{ass:comm} (Bounded Commutators).
        The commutators of the neighboring
        interaction and single site operators
        are densely defined and
        uniformly bounded\footnote{E.g., think of
        the canonical commutation relation for position
        and momentum operators.},
        i.e.,
        \begin{align*}
            \norm{\comm{\Phi_{j,j+1}}{H_{j+1}}}[\Lnorm]\leq J\quad\text{and}\quad
            \norm{\comm{H_j}{\Phi_{j, j+1}}}[\Lnorm]\leq J,
        \end{align*}
        for all $1\leq j\leq d-1$,
        where $\comm{A}{B}:=AB-BA$.

    \item\label{ass:self} (Self-Adjoint).
        The interaction and single site operators $\Phi_{j,j+1}$ and
        $H_{j}$ are self-adjoint.
    \end{enumerate}
\end{assumption}

\begin{remark}\label{rem:ass}
    Assumption \eqref{ass:local} means we only consider local 2-site
    interactions. Our results would remain unchanged for $k$-site interactions,
    for a fixed $k$. The point is that the complexity of approximating
    an eigenfunction scales exponentially with $k$ and not $d$. Moreover,
    we expect similar
    results could be obtained for long range interactions that decay
    sufficiently fast.

    We require Assumption \eqref{ass:self} since the proof heavily relies
    on the spectral decomposition. One could possibly
    extend the proofs presented
    here to sectorial operators. We are not certain to what extent
    approximability actually depends on the form
    of the resolvent/spectrum of the
    operator in $\C$.

    Assumption \eqref{ass:gap} is necessary for an area law to hold.
    Systems with degenerate ground states are at a
    quantum critical point and have been observed to exhibit
    divergent entanglement entropies (see
    \cite{PlenioAL, VidalQCP, Calabrese}).
    Finite-dimensionality of the GS eigenspace is required to ensure the GSP
    is trace-class-normalizable, otherwise the GSP makes no sense.

    Assumptions \eqref{ass:int} and \eqref{ass:comm} are required for the
    application of Lieb-Robinson bounds, i.e., finite speed information
    propagation. There are essentially two difficulties when considering
    information propagation for dynamics prescribed by an unbounded
    operator.

    First, unlike with classic
    Lieb-Robinson bounds (see \cite{LBOriginal}),
    bounded local operators do not have to remain local when evolved via the
    unitary operator $\exp{\i Ht}$ (see \cite{SuperSonice}).
    This can be remedied as in \cite{LBInfDim,
    NearComm} by,
    e.g., assuming the interactions in $H$ are of a certain type:
    bounded, as in this work, or specific
    types of unbounded operators that we do not consider here.
    Hence, we require Assumption \eqref{ass:int}.

    Second, when applying time dynamics to an unbounded local operator, it is
    not clear in which sense the operator remains \emph{approximately}
    local. Thus, Assumption \eqref{ass:comm} ensures that the non-local part
    is bounded.

    However, we essentially require only an application of
    Lieb-Robinson. Although Assumptions \eqref{ass:int} and \eqref{ass:comm}
    are certainly sufficient, they are perhaps not necessary.
\end{remark}

\begin{example}[Nearest Neighbor Interaction (NNI)]\label{ex:nni}
    We provide an example of how the general structure of an NNI
    Hamiltonian might look like.
    Perhaps the most famous example of an NNI Hamiltonian
    is the \emph{Ising model} (see \cite{Ising}).

    In this work we
    consider infinite-dimensional Hilbert spaces and
    unbounded Hamiltonians.
    A typical example to keep in mind is
    $\X=\bigotimes_{j=1}^d\X_j=\bigotimes_{j=1}^d\Lp{\R^n}{\C}$,
    where
    $n\in\set{1,\,2,\,3}$ if $H$ is to model a physical phenomenon.

    Let the Hamiltonian operator be given
    as
    \begin{align*}
        H=-\Delta+V.
    \end{align*}
    The Laplacian $\Delta$ is the one-site unbounded operator
    where $H_j=-\frac{1}{2}\frac{\partial^2}{\partial x_j^2}$. The potential
    $V$ contains the bounded interaction operators. E.g.,
    $V=\sum_{j=1}^{d-1}\Phi_{j,j+1}$, where $\Phi_{j,j+1}:\X\rightarrow\X$ is
    a bounded operator such as
    \begin{align*}
        (\Phi_{j,j+1}\uu)(x)&=c(x_j,\,x_{j+1})\uu(x),\quad\text{or}\\
        (\Phi_{j,j+1}\uu)(x)&=\int_{\R^{2n}}\kappa(x_j,\,x_{j+1},\,y_j,\,y_{j+1})
        \uu(x_1,\,\ldots,\,y_j,\,y_{j+1},\,\ldots,\,x_d)\d (y_j,\,y_{j+1}),
    \end{align*}
    where $c(\cdot)$ is a bounded coefficient function and $\kappa(\cdot)$ is an
    integral kernel.

    One would have to check the gap property and the degeneracy of the GS.
    Spectral properties and degeneracy of GSs
    have been extensively studied before and we refer to, e.g.,
    \cite[Chapter XIII]{MathPhys4} for more details.
\end{example}

We begin with a lemma that shows how we can approximately express the
GSP through the Hamiltonian operator. This will provide the
necessary link between the local operator structure and the local structure
of the GSP.

\begin{theoremEnd}{lemma}\label{lemma:pq}
    Let Assumption \ref{ass:op} \eqref{ass:gap} hold.
    Then, for any $q>0$ and
    \begin{align}\label{eq:pq}
        \rho^q:=\frac{1}{N\sqrt{2\pi q}}\int_{-\infty}^\infty
        \exp{\i Ht}\exp{-\frac{t^2}{2q}}\d t,
    \end{align}
    we have
    \begin{align*}
        \norm{\rho^q-P_0}[\Lnorm]
        \leq\exp{-\frac{1}{2}(\Delta E)^2q},
    \end{align*}
    with $\Delta E$ from \eqref{eq:spgap}.
\end{theoremEnd}
\begin{proofEnd}
    The operator $U(t):=\exp{\i Ht}\exp{-\frac{t^2}{2q}}$ is strongly continuous
    for all $t\in\R$. Thus, a finite integral of $U(t)$ is well defined.
    For any $\uu\in\X$
    \begin{align*}
        \lim_{c\rightarrow\infty}\norm{\frac{1}{\sqrt{2\pi q}}
        \int_{-c}^cU(t)\uu\d t}[\X]\leq
        \lim_{c\rightarrow\infty}\norm{\uu}[\X]\frac{1}{\sqrt{2\pi q}}
        \int_{-c}^c\exp{-\frac{t^2}{2q}}\d t=\norm{\uu}[\X].
    \end{align*}
    Thus, the integral \eqref{eq:pq} is well defined.

    Since $H$ is self-adjoint, we have the spectral decomposition
    \begin{align*}
        H=\int_{\sigma(H)}\lambda\d P(\lambda),
    \end{align*}
    where $P:\sigma(H)\rightarrow\LL{\X}$ is a projection valued measure.
    Due to the gap assumption, we get that $NP_0=P(\lambda_0)$.

    Applying functional calculus for self-adjoint operators
    \begin{align*}
        \exp{\i Ht}=\int_{\sigma(H)}\exp{\i\lambda t}\d P(\lambda).
    \end{align*}
    Equation
    \eqref{eq:pq} is to be interpreted as the unique operator such that
    for any $\uu\in\X$
    \begin{align*}
        \inp{\uu}{\frac{1}{\sqrt{2\pi q}}\int_{-\infty}^\infty U(t)\uu\d t}_{\X}
        &=
        \frac{1}{\sqrt{2\pi q}}\int_{-\infty}^\infty
        \exp{-\frac{t^2}{2q}}\inp{\uu}{U(t)\uu}_{\X}\d t\\
        &=
        \frac{1}{\sqrt{2\pi q}}\int_{-\infty}^\infty
        \exp{-\frac{t^2}{2q}}\int_{\sigma(H)}\exp{\i\lambda t}\d P_{\uu}(\lambda)\d t,
    \end{align*}
    where $P_{\uu}(\cdot)=\inp{\uu}{P(\cdot)\uu}_{\X}$ and the equality
    follows from the linearity and continuity of the $\X$-inner product.
    For the last integral we can apply Fubini's Theorem
    for general product measures.
    This allows us to write
    \begin{align*}
        \rho^q&=\frac{1}{N\sqrt{2\pi q}}\int_{-\infty}^\infty \int_{\sigma(H)}
        \exp{\i\lambda t}\exp{-\frac{t^2}{2q}}\d P(\lambda)\d t\\
        &\overset{\text{gap}}{=}\frac{1}{N\sqrt{2\pi q}}\int_{-\infty}^\infty
        P_0\exp{-\frac{t^2}{2q}}+\int_{\sigma(H)\setminus\set{\lambda_0}}
        \exp{i\lambda t}\exp{-\frac{t^2}{2q}}\d P(\lambda)\d t\\
        &\overset{\text{Fubini}}{=}
        P_0+\frac{1}{N\sqrt{2\pi q}}\int_{\sigma(H)\setminus\set{\lambda_0}}
        \int_{-\infty}^\infty\exp{\i\lambda t}\exp{-\frac{t^2}{2q}}\d t\d P(\lambda).
    \end{align*}
    The last term is the Fourier transform of the density of the
    normal distribution. Thus,
    \begin{align*}
        \norm{\rho^q-P_0}[\Lnorm]=\frac{1}{N}
        \norm{\int_{\sigma(H)\setminus\set{\lambda_0}}
        \exp{-\frac{1}{2}\lambda^2 q}\d P(\lambda)}[\Lnorm]
        \leq\frac{1}{N}\exp{-\frac{1}{2}(\Delta E)^2q},
    \end{align*}
    which completes the proof.
\end{proofEnd}

Next, we want to approximate $H$ by a sum of three local operators, where each
operator
approximately annihilates the ground state. To this end, we apply Hasting's
quasi-adiabatic continuation
technique (see \cite{QA, Hastings}),
which was also studied
in \cite{NearComm} in the infinite-dimensional setting\footnote{There the
authors considered this technique in order to describe states that belong to
the same phase.}.

\begin{theoremEnd}{lemma}\label{lemma:ann}
    Suppose Assumption \ref{ass:op} (\eqref{ass:local}-\eqref{ass:comm})
    holds.

    For a fixed $l\in\N$ and a
    fixed $1+l\leq j\leq d-2-l$,
    \begin{align*}
        H_L&:=\sum_{k\leq j-l-2}H_{k,k+1},\quad
        H_B:=\sum_{j-l-1\leq k\leq j+l+1}H_{k,k+1},\quad
        H_R:=\sum_{k\geq j+l+2}H_{k,k+1}.
    \end{align*}
    Let
    $\bs\psi_0$ be an arbitrary GS and    
    w.l.o.g.\footnote{Otherwise the integrals can be
    re-centered around the expectation value.} let $\inp{\uu_0}{H_L\uu_0}_{\X}=
    \inp{\uu_0}{H_B\uu_0}_{\X}=\inp{\uu_0}{H_R\uu_0}_{\X}=0$.
    Then, for any $q>0$ and
    \begin{align*}
        \tilde{H}_L&:=\frac{1}{\sqrt{2\pi q}}\int_{-\infty}^\infty H_L(t)\exp{-\frac{t^2}{2q}}\d t,\\
        \tilde{H}_B&:=\frac{1}{\sqrt{2\pi q}}\int_{-\infty}^\infty H_B(t)\exp{-\frac{t^2}{2q}}\d t,\\
        \tilde{H}_R&:=\frac{1}{\sqrt{2\pi q}}\int_{-\infty}^\infty H_R(t)\exp{-\frac{t^2}{2q}}\d t,
    \end{align*}
    where $H_{\cdots}(t)$ is
    the Heisenberg evolution of $H_{\cdots}$ defined as
    \begin{align*}
        H_L(t):=\exp{\i Ht}H_L\exp{-\i Ht},
    \end{align*}
    we have
    \begin{align*}
        \norm{\tilde{H}_L\uu_0}[\X]&\leq
        3J^2(\Delta E)^{-1}\exp{-\frac{1}{2}(\Delta E)^2q},\\
        \norm{\tilde{H}_B\uu_0}[\X]&\leq
        3J^2(\Delta E)^{-1}\exp{-\frac{1}{2}(\Delta E)^2q},\\
        \norm{\tilde{H}_R\uu_0}[\X]&\leq
        3J^2(\Delta E)^{-1}\exp{-\frac{1}{2}(\Delta E)^2q}.
    \end{align*}
    The constant $J$ is the interaction strength from \eqref{eq:intstr}.
\end{theoremEnd}
\begin{proofEnd}
    By the same arguments as in Lemma \ref{lemma:pq}, the integrals are
    well defined. Next, application to the ground state yields
    \begin{align*}
        H_L(t)\uu_0&=\exp{\i Ht}H_L\int_{\sigma(H)}\exp{-\i\lambda t}\d P(\lambda)
        \uu_0\overset{\lambda_0=0}{=}\exp{\i Ht}H_L\uu_0
        =
        \int_{\sigma(H)}\exp{\i\lambda t}\d P(\lambda)H_L\uu_0\\
        &=\inp{\uu_0}{H_L\uu_0}_{\H}\uu_0+\int_{\sigma(H)\setminus\set{\lambda_0}}
        \exp{\i\lambda t}\d P(\lambda)H_L\uu_0.
    \end{align*}
    Thus,
    \begin{align*}
        \tilde{H}_L\uu_0&=\frac{1}{\sqrt{2\pi q}}\int_{-\infty}^\infty
        H_L(t)\uu_0\exp{-\frac{t^2}{2q}}\d t\\
        &=\int_{\sigma(H)\setminus\set{\lambda_0}}\frac{1}{\sqrt{2\pi q}}
        \int_{-\infty}^\infty\exp{\i\lambda t}\exp{-\frac{t^2}{2q}}\d t\d P(\lambda)
        H_L\uu_0\\
        &=
        \int_{\sigma(H)\setminus\set{\lambda_0}}\exp{-\frac{1}{2}\lambda^2 q}
        \d P(\lambda)H_L\uu_0.
    \end{align*}
    On the other hand,
    \begin{align*}
        H\tilde{H}_L\uu_0=\int_{\sigma(H)\setminus\set{\lambda_0}}
        \lambda\exp{-\frac{1}{2}\lambda^2q}\d P(\lambda)H_L\uu_0.
    \end{align*}
    Hence,
    \begin{align}\label{eq:hhl}
        \norm{H\tilde{H}_L\uu_0}[\X]\geq \Delta E\norm{
        \int_{\sigma(H)\setminus\set{\lambda_0}}\exp{-\frac{1}{2}\lambda^2q}
        \d P(\lambda)H_L\uu_0}[\X]=\Delta E\norm{\tilde{H}_L\uu_0}[\X].
    \end{align}
    Next, since $H\uu_0=0$,
    \begin{alignat*}{2}
        H\tilde{H}_L\uu_0&=&&(H\tilde{H}_L-\tilde{H}_LH)\uu_0
        =\int_{\sigma(H)\setminus\set{\lambda_0}}\exp{-\frac{1}{2}\lambda^2q}
        \d P(\lambda)\comm{H}{H_L}\uu_0\\
        &=&&\int_{\sigma(H)\setminus\set{\lambda_0}}\exp{-\frac{1}{2}\lambda^2q}
        \d P(\lambda)\comm{H_{j-l-1,j-l}}{H_{j-l-2,j-l-1}}\uu_0\\
        &=&&\int_{\sigma(H)\setminus\set{\lambda_0}}\exp{-\frac{1}{2}\lambda^2q}
        \d P(\lambda)
        \bigg\{
        \comm{H_{j-l-1}}{(\Phi_{j-l-2,j-l-1}-\Phi_{j-l-1,j-l})}\\
        & &&+\comm{\Phi_{j-l-1,j-l}}{\Phi_{j-l-2,j-l-1}}\bigg\}\uu_0,
    \end{alignat*}
    where the last equality follows from Assumption \eqref{ass:int}.
    And thus
    \begin{align*}
        \norm{H\tilde{H}_L\uu_0}[\X]\leq 3J^2\exp{-\frac{1}{2}(\Delta E)^2q}.
    \end{align*}
    Together with \eqref{eq:hhl}
    \begin{align*}
        \norm{\tilde{H}_L\uu_0}[\X]\leq
        3J^2(\Delta E)^{-1}\exp{-\frac{1}{2}(\Delta E)^2q},
    \end{align*}
    and analogously for $\tilde{H}_B$, $\tilde{H}_R$. This completes the
    proof.
\end{proofEnd}

\begin{remark}
    Note that the above bound depends explicitly only on $q$, $\Delta E$
    and $J$. In fact, more
    precisely, the bound depends on an estimate for
    \begin{align*}
        \norm{\comm{H}{H_L}}[\Lnorm],
    \end{align*}
    and the latter was assumed in \eqref{eq:intstr} to be uniformly bounded,
    i.e., in particular it is independent of $j$ or $l$.

    However, the subsequent lemmas will employ Lieb-Robinson bounds that
    depend explicitly on the parameter $l$. Thus, we will eventually use the
    above lemma and choose
    the constant $q$ depending on the spectral gap
    $\Delta E$ and the parameter $l$.
\end{remark}

Next, we show that $\tilde{H}_L$, $\tilde{H}_B$ and $\tilde{H}_R$
are approximately local\footnote{In the sense specified by the following lemma.}.
This is mainly due to
Lieb-Robinson type estimates.
\begin{theoremEnd}{lemma}\label{lemma:ml}
    Under Assumption \ref{ass:op} (\eqref{ass:local}, \eqref{ass:int},
    \eqref{ass:comm}), there exist local bounded operators
    $\diff_L$, $\diff_B$ and
    $\diff_R$ supported on
    $\X_{j-2l-2, j}$,
    $\X_{j-2l-2,j+2l+3}$ and
    $\X_{j+1, j+2l+3}$,
    respectively, such that for
    \begin{align*}
        M_L&:=H_L+\diff_L,\quad
        M_B:=H_B+\diff_B,\quad
        M_R:=H_R+\diff_R,
    \end{align*}
    there exist constants $c_1>0$, $C_1>0$ such that
    \begin{align*}
        \norm{\tilde{H}_L-M_L}[\Lnorm]&\leq C_1J^2\max\set{q^{1/2},\,q^{3/2}}\exp{-c_1l},\\
        \norm{\tilde{H}_B-M_B}[\Lnorm]&\leq C_1J^2\max\set{q^{1/2},\,q^{3/2}}\exp{-c_1l},\\
        \norm{\tilde{H}_R-M_R}[\Lnorm]&\leq C_1J^2\max\set{q^{1/2},\,q^{3/2}}\exp{-c_1l},
    \end{align*}
    where $q$, $l$ are the parameters from Lemma \ref{lemma:ann} and
    $J$ is the interaction strength from \eqref{eq:intstr}.
\end{theoremEnd}
\begin{proofEnd}
    First, note that we can differentiate $H_L(t)$ to obtain
    \begin{align*}
        \frac{\d}{\d t}H_L(t)&=\frac{\d}{\d t}\exp{\i Ht}H_L\exp{-\i Ht}\\
        &=\exp{\i Ht}\i HH_L\exp{-\i Ht}-\exp{\i Ht}iH_LH\exp{-\i Ht}\\
        &=\exp{\i Ht}\i\comm{H}{H_L}\exp{-\i Ht}
        =:\i\comm{H}{H_L}(t),
    \end{align*}
    and $H_L(0)=H_L$.
    Thus, we can write
    \begin{align}\label{eq:thanks1}
        H_L(t)=H_L+\int_{0}^t\i\comm{H}{H_L}(\tau)\d\tau.
    \end{align}
    By Assumption, $\comm{H}{H_L}$ is bounded and supported on
    $\X_{j-l-2,j-l}$. Consequently, we can write
    \begin{align*}
        \tilde{H}_L=H_L+\frac{1}{\sqrt{2\pi q}}\int_{-\infty}^\infty
        \int_0^t\i\comm{H}{H_L}(\tau)\d\tau\exp{-\frac{t^2}{2q}}\d t.
    \end{align*}
    Since the commutator is bounded and local, and the interactions in
    $H$ are bounded, by \cite[Corollary 2.2]{LBInfDim},
    we know a Lieb-Robinson bound
    applies to $\comm{H}{H_L}$. I.e., there exists a constant (velocity)
    $v\geq 0$ and constants $C>0$, $a>0$, such that
    \begin{align}\label{eq:lrbnd}
        &\norm{\comm{\comm{H}{H_L}(\tau)}{B}}[\Lnorm]
        \leq C\norm{\comm{H}{H_L}}[\Lnorm]
        \norm{B}[\Lnorm]\exp{-a\left\{\dist(\comm{H}{H_L}, B)-v|\tau|
        \right\}},
    \end{align}
    for all bounded and local $B$.
    Thus, by \cite[Lemma 3.2]{NearComm},
    there exists a map $\Pi:\LL{\X}\rightarrow\LL{\X}$ such that
    $\Pi(A)$ is supported on $\X_{j-2l-2,j}$ and, for any $A\in\LL{\X}$
    satisfying \eqref{eq:lrbnd} with $\dist(A, B)\geq l$, we have
    \begin{align*}
        \norm{A-\Pi(A)}[\Lnorm]\leq 2C\norm{A}[\Lnorm]
        \exp{-a\left\{l-v|\tau|
        \right\}}.
    \end{align*}
    Then, using \cite[Theorem 3.4]{NearComm}, we integrate over time
    and further estimate
    \begin{align*}
        &\norm{\int_{0}^t\comm{H}{H_L}(\tau)\d\tau-
        \int_0^t\Pi\left(\comm{H}{H_L}(\tau)\right)
        \d\tau}[\Lnorm]
        \leq
        |t|CJ^2\exp{-a\{l-v|t|\}}.
    \end{align*}

    We define the operator
    \begin{align*}
        \diff_L:=\frac{1}{\sqrt{2}\pi q}\int_{-\infty}^\infty
        \int_0^t \Pi\left(\i\comm{H}{H_L}(\tau)\right)\d\tau\exp{-\frac{t^2}{2q}}
        \d t.
    \end{align*}
    By all of the above,
    this operator is bounded and supported on $\X_{j-2l-2, j}$. Analogously,
    we define $\diff_B$ and $\diff_R$ with supports in
    $\X_{j-2l-2,j+2l+3}$ and $\X_{j+1, j+2l+3}$, respectively.

    What remains is to truncate the tails of the integral to obtain an
    overall error of the same order as the Lieb-Robinson bound. Let
    $T=\frac{l}{2v}$. Then,
    \begin{align}\label{eq:trunc}
        &\norm{\tilde{H}_L^1-M_L}[\Lnorm]=\notag\\
        &=
        \norm{\frac{1}{\sqrt{2\pi q}}\int_{-\infty}^\infty\int_{0}^t
        \left(\i\comm{H}{H_L}(\tau)-
        \Pi\left(\i\comm{H}{H_L}(\tau)\right)\right)\d\tau
        \exp{-\frac{t^2}{2q}}\d t}[\Lnorm]\notag\\
        &=\norm{\frac{1}{\sqrt{2\pi q}}\left(
        \int_{|t|\leq T}\ldots+\int_{|t|>T}\ldots\right)}[\Lnorm]\notag\\
        &\leq
        CJ^2\exp{-al}
        \left(\exp{avT}\frac{1}{\sqrt{2\pi q}}\int_{|t|\leq T}
        |t|\exp{-\frac{t^2}{2q}}\d t+\frac{1}{\sqrt{2\pi q}}\int_{|t|>T}
        |t|\exp{av|t|-\frac{t^2}{2q}}\d t\right).
    \end{align}
    For the first term
    \begin{align*}
        \frac{1}{\sqrt{2\pi q}}\int_{|t|\leq T}
        |t|\exp{-\frac{t^2}{2q}}\d t\leq \sqrt{\frac{q}{2\pi}}.
    \end{align*}
    For the second
    \begin{alignat*}{2}
        \frac{1}{\sqrt{2\pi q}}\int_{|t|>T}
        |t|\exp{av|t|-\frac{t^2}{2q}}\d t
        &=&&\sqrt{\frac{2}{\pi q}}\int_{T}^\infty
        t\exp{avt-\frac{t^2}{2q}}\d t\\
        &=&&q\exp{avT-\frac{T^2}{2q}}\\
        & &&+avq
        \int_T^\infty\exp{avt-\frac{t^2}{2q}}\d t.
    \end{alignat*}
    For the latter term
    \begin{align*}
        \int_T^\infty\exp{avt-\frac{t^2}{2q}}\d t&=
        \int_T^\infty\exp{\frac{(av)^2q}{2}-
        \left(\sqrt{\frac{1}{2q}}t-av\sqrt{\frac{q}{2}}\right)^2}\d t\\
        &=
        \sqrt{2q}\int_{\sqrt{\frac{1}{2q}}T-av\sqrt{\frac{q}{2}}}^\infty
        \exp{-\tau^2}\d\tau\leq
        \sqrt{\frac{q\pi}{2}}\exp{avT-\frac{T^2}{2q}}.
    \end{align*}
    And thus
    \begin{alignat*}{2}
        \norm{\tilde{H}_L-M_L}[\Lnorm]
        &\leq&&
        CJ^2\exp{-al}\times\\
        & &&\bigg(
        \exp{al/2}\sqrt{\frac{q}{2\pi}}
        +q\exp{al/2}\exp{-\frac{T^2}{2q}}
        +avq^{3/2}\sqrt{\frac{\pi}{2}}\exp{al/2}\exp{-\frac{T^2}{2q}}
        \bigg)\\
        &\leq&& CJ^2\exp{-\frac{a}{2}l}
        \max\set{
        \sqrt{\frac{q}{2\pi}},\, q,\, avq^{3/2}\sqrt{\frac{\pi}{2}}}
        \left(1+2\exp{-\frac{l^2}{8v^2q}}\right)\\
        &\leq&& C_1J^2\max\set{q^{1/2},\,q^{3/2}}\exp{-c_1l},
    \end{alignat*}
    with $C_1$ and $c_1$ defined in an obvious way as above.
    This completes the proof.
\end{proofEnd}

We can conclude the existence of the first two operators that
we will need to approximate $P_0$.

\begin{theoremEnd}{lemma}\label{lemma:ol}
    Under Assumption \ref{ass:op}, there exist local, bounded and
    self-adjoint (projection) operators $L=L(j,l)$,
    $R=R(j,l)$ with the property
    \begin{align*}
        \norm{(L-\id)\uu_0}[\X]&\leq \exp{-c_1l/2},\\
        \norm{(R-\id)\uu_0}[\X]&\leq \exp{-c_1l/2},
    \end{align*}
    for any $GS$ $\bs\psi_0$.
    The operators $L$ and $R$ have the same support as $M_L$
    and $M_R$, respectively, and $\norm{L}[\Lnorm]
    =\norm{R}[\Lnorm]=1$.
\end{theoremEnd}
\begin{proofEnd}
    Recall from Lemma \ref{lemma:ml} that we applied Lieb-Robinson to operators
    such as $\i\comm{H}{H_L}(\tau)$. By Assumption \ref{ass:op},
    since $H$ and $H_L$ are self-adjoint, the commutators are bounded and
    one has
    $\adj{\i\comm{H}{H_L}}=-\i\left(H_LH-HH_L\right)=\i\comm{H}{H_L}$
    $\Rightarrow$ we
    can conclude
    that the commutator is self-adjoint. Applying Lieb-Robinson as in
    \cite[Lemma 3.2]{NearComm}, one can construct the local approximation to
    preserve self-adjointness. We briefly elaborate.

    Let $\X=\X_1\otimes \X_2$.
    For the construction of $\Pi:\LL{\X}\rightarrow\LL{\X}$
    in \cite[Lemma 3.2]{NearComm}, the authors use an arbitrary
    state\footnote{Positive operator with unit trace norm.}
    $\rho:\X_2\rightarrow\X_2$, though the resulting bound does not depend on the
    choice of $\rho$. Then, for this state, applying the spectral decomposition,
    we have
    \begin{align*}
        \rho=\sum_{k=1}^\infty\lambda_k\inp{\cdot}{\bs\psi_k}_{\X_2}\bs\psi_k,
    \end{align*}
    for $\set{\bs\psi_k:\;k\in\N}$ orthonormal in $\X_2$. Next, for each $\bs\psi_k$,
    the authors define the map $A_k$ as
    \begin{align*}
        \inp{\bs v}{A_k\bs w}_{\X_1}=\inp{\bs v\otimes\bs\psi_k}{A\bs w\otimes\bs\psi_k}_{\X},\quad
        \forall \bs v,\,\bs w\in\X_1.
    \end{align*}
    Finally, the map $\Pi(A)$ is defined as
    \begin{align*}
        \Pi(A):=\left(\sum_{k=1}^\infty\lambda_kA_k \right)\otimes\id.
    \end{align*}
    Note that each $A_k$ is self-adjoint if $A$ is self-adjoint. Therefore,
    $\Pi(A)$ is self-adjoint.

    Thus, $M_L$, $M_B$ and $M_R$ can be chosen
    self-adjoint.
    By Lemmas \ref{lemma:ml} and \ref{lemma:ann}, picking
    $q=c_1\frac{2l}{(\Delta E)^2}$, we get
    $\norm{M_L\uu_0}[\X]\leq C_1J^2\max\set{q^{1/2},\,q^{3/2}}\exp{-c_1l}$.
    Moreover, since $M_L$ is
    self-adjoint, there exists a projection valued measure $P(\cdot)$ such that
    \begin{align*}
        \normsq{M_L\uu_0}{\X}=\inp{M_L\uu_0}{M_L\uu_0}_{\X}=
        \inp{\uu_0}{M_L^2\uu_0}_{\X}=\int_{\sigma(M_L)}\lambda^2\d
        P_{\uu_0}(\lambda).
    \end{align*}

    We split the spectrum of $M_L$ as
    \begin{align*}
        \sigma_1(M_L)&:=\set{\lambda\in\sigma(M_L):|\lambda|\leq
        C_1J^2\max\set{q^{1/2},\,q^{3/2}}\exp{-c_1l/2}},\\
        \sigma_2(M_L)&:=\sigma(M_L)\setminus
        \sigma_1(M_L).
    \end{align*}
    Define $L$ as
    \begin{align*}
        L:=\int_{\sigma_1(M_L)}\d P(\lambda).
    \end{align*}
    Clearly, $L$ is a bounded self-adjoint operator with
    $\norm{L}[\Lnorm]=1$ and the same support as $M_L$.
    Moreover, by orthogonality of the spectral subspaces
    \begin{align*}
        C_1^2J^4\max\set{q,\,q^{3}}\exp{-c_12l}&\geq \norm{M_L\uu_0}[\X]
        =
        \int_{\sigma_1(M_L)}\lambda^2\d P_{\uu_0}(\lambda)\uu_0+
        \int_{\sigma_2(M_L)}\lambda^2\d P_{\uu_0}(\lambda)\uu_0\\
        &\geq C_1^2 J^4\max\set{q,\,q^{3}}
        \exp{-c_1l}\int_{\sigma_2(M_L)}\d P_{\uu_0}(\lambda).
    \end{align*}
    Thus, $\int_{\sigma_2(M_L)}\d P_{\uu_0}(\lambda)\leq\exp{-c_1l}$.
    Finally, this gives
    \begin{align*}
        \norm{(L-\id)\uu_0}[\X]&=
        \norm{\int_{\sigma_2(M_L)}\d P(\lambda)\uu_0}[\X]=
        \left(\int_{\sigma_2(M_L)}\d P_{\uu_0}(\lambda)\right)^{1/2}\\
        &\leq\exp{-c_1l/2}.
    \end{align*}
    Analogously for $R$. This completes the proof.
\end{proofEnd}

What remains is a step by step approximation of $\rho^q$ as a product of three
local operators.

\begin{theoremEnd}{lemma}\label{lemma:tpq}
    Under Assumption \ref{ass:op}, we can further approximate $P_0$ as
    \begin{align*}
        \tilde{\rho}^q&:=\frac{1}{N\sqrt{2\pi q}}\int_{-\infty}^\infty
        \adj{\oe{A}{0}{t}}\exp{-\frac{t^2}{2q}}LR\d t,\\
        A(t)&:=\exp{\i (M_L+M_R)t}\i M_B\exp{-\i (M_L+M_R)t},
    \end{align*}
    where $\adj{\oe{A}{0}{t}}$ is the negative time-ordered
    exponential
    and $q=c_1\frac{2l}{(\Delta E)^2}$. We have
    \begin{align*}
        \norm{\tilde{\rho}^q-P_0}[\Lnorm]\leq
        C_2(1/N)J^2\max\set{q,\,q^2}\left(
        \exp{-c_1l/2}+\exp{-c_1l}\right),
    \end{align*}
    for some constant $C_2>0$.
\end{theoremEnd}
\begin{proofEnd}
    Note that $H=\tilde{H}_L+\tilde{H}_B+
    \tilde{H}_R$. By utilizing the estimates from Lemma \ref{lemma:ml}, we have
    \begin{align}\label{eq:exp}
        &\norm{\frac{1}{N\sqrt{2\pi q}}\int_{-\infty}^\infty
        \big\{\exp{\i(\tilde{H}_L+\tilde{H}_B+\tilde{H}_R
        )t}\notag
        -\exp{\i(M_L+M_B+M_R)t}\big\}\exp{-\frac{t^2}{2q}}\d t}[\Lnorm]\notag\\
        &\leq 3C_1(1/N)J^2\max\set{q^{1/2},\, q^{3/2}}\exp{-c_1l}
        \frac{1}{\sqrt{2\pi q}}\int_{-\infty}^\infty |t|
        \exp{-\frac{t^2}{2q}}\d t\\
        &=3(2\pi)^{-1/2}C_1(1/N)J^2
        \max\set{q,\, q^{2}}\exp{-c_1l}\notag,
    \end{align}
    where the inequality can be shown using
    \cite[Chapter 9, Theorem 2.12, Equation (2.22)]{Kato}.

    Next, for the exponential term we can write
    \begin{align*}
        &\exp{\i(M_L+M_B+M_R)t}
        =\exp{\i(M_L+M_B+M_R)t}\exp{-\i(M_L+M_R)t}
        \exp{\i(M_L+M_R)t},
    \end{align*}
    and we define $U(t):=\exp{\i(M_L+M_B+M_R)t}\exp{-\i(M_L+M_R)t}$.
    In $U(t)$ the term in the exponent commutes for different $t$, since
    the time-dependence is a simple multiplication by $t$. We thus
    compute
    \begin{align*}
        &\frac{\d}{\d t}U(t)
        =\exp{\i(M_L+M_B+M_R)t}\i(M_L+M_B+M_R)
        \exp{-\i(M_L+M_R)t}\\
        &-
        \exp{\i(M_L+M_B+M_R)t}\i(M_L+M_R)
        \exp{-\i(M_L+M_R)t}\\
        &=\exp{\i(M_L+M_B+M_R)t}\i M_B\exp{-\i(M_L+M_R)t}\\
        &=\exp{\i(M_L+M_B+M_R)t}\exp{-\i(M_L+M_R)t}
        \exp{\i(M_L+M_R)t}\i M_B\exp{-\i(M_L+M_R)t}\\
        &=U(t)\exp{\i(M_L+M_R)t}\i M_B\exp{-\i(M_L+M_R)t}
    \end{align*}
    and $U(0)=\id$. We abbreviate
    \begin{align*}
        M_B(t):=\exp{\i(M_L+M_R)t}M_B\exp{-\i(M_L+M_R)t}.
    \end{align*}
    Due to the simple form of $\i M_B(t)$, the solution to this initial
    value problem exists, is unique and is
    given by the (negative) time-ordered exponential of
    $\i M_B(t)$, see \cite[Chapter X.12]{MathPhys2}
    (interaction representation) or \cite{Unbounded} for
    ordered exponentials of more general unbounded time-dependent
    Hamiltonians.

    Thus, our approximation so far is
    \begin{align*}
        \bar{\rho}^q=\frac{1}{N\sqrt{2\pi q}}\int_{-\infty}^\infty
        \adj{\oe{A}{0}{t}}\exp{\i(M_L+M_R)t}\exp{-\frac{t^2}{2q}}\d t.
    \end{align*}
    By multiplying $L$, $R$ from the right we obtain
    \begin{align*}
        &\norm{\bar{\rho}^qLR-P_0}[\Lnorm]
        =\norm{(\bar{\rho}^q-P_0+P_0)LR-P_0}[\Lnorm]
        =\norm{(\bar{\rho}^q-P_0)LR+P_0LR-P_0}[\Lnorm]
        \\
        &\leq
        \norm{\bar{\rho}^q-P_0}[\Lnorm]+
        \norm{P_0LR-P_0}[\Lnorm].
    \end{align*}
    For the latter term we use Lemma \ref{lemma:ol}, set
    $q=c_1\frac{2l}{(\Delta E)^2}$ and obtain
    \begin{align*}
        &\norm{P_0LR-P_0}[\Lnorm]=
        \norm{P_0
        \left\{(L-\id+\id)+(R-\id+\id)-\id\right\}
        }[\Lnorm]\\
        &=\norm{P_0(L-\id)(R-\id)+P_0(L-\id)+
        P_0(R-\id)
        }[\Lnorm]\\
        &\leq
        3\norm{(L-\id)P_0}[\Lnorm]+
        \norm{(R-\id)P_0}[\Lnorm]
        \leq 4\exp{-c_1l/2}.
    \end{align*}
    since
    $\norm{P_0}[\Lnorm]=
    \norm{L}[\Lnorm]=1$ and $P_0$, $L$ and $R$ are self-adjoint.
    Thus, overall
    \begin{align*}
        \norm{\bar{\rho}^qLR-P_0}[\Lnorm]&\leq
        3(2\pi)^{-1/2}(1/N)C_1J^2\max\set{q,\, q^2}\exp{-c_1l}
        +4\exp{-c_1l/2}.
    \end{align*}

    Finally, by definition, $L$ and $R$ project onto the spectral subspaces
    of $M_L$ and $M_R$ corresponding to small eigenvalues. I.e.,
    \begin{align*}
        \norm{\left[\exp{\i(M_L+M_R)t}-\id\right]LR
        }[\Lnorm]\leq 2|t|C_1J^2\max\set{q^{1/2},\,q^{3/2}}
        \exp{-c_1l/2},
    \end{align*}

    hence
    \begin{align*}
        \norm{\bar{\rho}^qLR-\tilde{\rho}^q}[\Lnorm]
        &\leq\bigg\|\frac{1}{N\sqrt{2\pi q}}\int_{-\infty}^\infty
        \adj{\oe{A}{0}{t}}\times\\
        &\exp{-\frac{t^2}{2q}}\left(\exp{\i(M_L+M_R)t}LR-LR
        \right)\d t\bigg\|_{\Lnorm}\\
        &\leq 2C_1(1/N)J^2\max\set{q^{1/2},\, q^{3/2}}\exp{-c_1l/2}
        \frac{1}{\sqrt{2\pi q}}\int_{-\infty}^\infty |t|\exp{-\frac{t^2}{2q}}\d t\\
        &=
        C_1(1/N)J^2(2\pi)^{-1/2}\max\set{q,\,q^2}\exp{-c_1l/2}.
    \end{align*}

    Overall we obtain
    \begin{align*}
        \norm{\tilde{\rho}^q-P_0}[\Lnorm]&\leq
        \norm{\bar{\rho}^qLR-P_0}[\Lnorm]+
        \norm{\bar{\rho}^qLR-\tilde{\rho}^q}[\Lnorm]\\
        &\leq C_2J^2\max\set{q,\, q^2}\left(
        \exp{-c_1l/2}+\exp{-c_1l}\right),
    \end{align*}
    for an appropriately chosen constant $C_2>0$. This completes the proof.
\end{proofEnd}

It remains to show how we can obtain a local operator $B$, maintaining the
same approximation order. This follows once more from a Lieb-Robinson bound.

\begin{theoremEnd}{lemma}\label{lemma:ob}
    Consider the operator
    \begin{align*}
        \OB:=\frac{1}{N\sqrt{2\pi q}}\int_{-\infty}^\infty
        \adj{\oe{A}{0}{t}}\exp{-\frac{t^2}{2q}}\d t,
    \end{align*}
    with $A(t)$ as above
    \begin{align*}
        A(t):=\exp{\i(M_L+M_R)t}\i M_B\exp{-\i(M_L+M_R)t}.
    \end{align*}
    Then, there exists a local bounded operator $B$
    supported on $\X_{j-3l-2,j+3l+3}$, with $\norm{B}[\Lnorm]\leq 1$
    such that
    \begin{align*}
        \norm{\OB-B}[\Lnorm]\leq
        C_3J^2\max\set{q^{1/2},\, q}\exp{-c_3l},
    \end{align*}
    for some constants $C_3>0$, $c_3>0$.
\end{theoremEnd}
\begin{proofEnd}
    We use the same technique as in Lemma \ref{lemma:ml} to show that the non-local
    part of $M_B(t)$ is bounded.
    \begin{align*}
        M_B(t)&:=\exp{\i(M_L+M_R)t}M_B\exp{\i(M_L+M_R)t},\\
        \frac{\d}{\d t}M_B(t)&=\exp{\i(M_L+M_R)t}\comm{M_L+M_R}{M_B}
        \exp{-\i(M_L+M_R)t},\\
        \comm{M_L+M_R}{M_B}&=\comm{H_L}{H_B}+\comm{H_L}{\diff_B}+
        \comm{\diff_L}{H_B}+\comm{\diff_L}{\diff_B}+
        \comm{H_R}{H_B}\\
        &+\comm{H_R}{\diff_B}+\comm{\diff_R}{H_B}
        +\comm{\diff_R}{\diff_B}.
    \end{align*}
    Since $\diff_L$ and $\diff_R$ are supported on a superset of the supports
    of $\comm{H}{H_L}$ and $\comm{H}{H_R}$, we obtain
    \begin{align*}
        \supp{\comm{M_L+M_R}{M_B}}
        \subset\supp{\diff_L}\cup\supp{\diff_R}=\supp{M_B}.
    \end{align*}
    Thus, as in Lemma \ref{lemma:ml},
    we can approximate $M_B(t)$ by
    a local operator $\tilde{M}_B(t)$ supported on $\X_{j-3l-2,j+3l+3}$
    such that
    \begin{align*}
        \norm{M_B(t)-\tilde{M}_B(t)}[\Lnorm]\leq
        CJ^2|t|\exp{-a\{l-v|t|\}}.
    \end{align*}
    We utilize this estimate in a similar way as in \eqref{eq:exp}. To this end,
    we write the time-ordered exponential as a product integral
    and use a
    step function approximation to the integral
    (see \cite[Chapter 3.6]{ProductInt}), namely
    \begin{align*}
        \adj{\oe{A}{0}{t}}=\prod_{0}^t\exp{A(\tau)}\d\tau=
        \lim_{K\rightarrow\infty}\left(
        \exp{A(t_K)\Delta t}\cdots\exp{A(t_0)\Delta t}\right),
    \end{align*}
    where $t_i=i\Delta t$, $\Delta t=t/K$ and the convergence is meant in
    the strong sense. This is possible due to the simple form of $A(t)$,
    i.e., $\exp{A(t)}$ is bounded with norm 1.

    For $K=1$, we obtain as in \eqref{eq:exp}
    \begin{align*}
        \norm{\exp{\i M_B(t)\Delta t}-\exp{\i\tilde{M}_B(t)\Delta t}
        }[\Lnorm]\leq |\Delta t|\norm{M_B(t)-\tilde{M}_B(t)
        }[\Lnorm],
    \end{align*}
    with $|\Delta t|=|t|$. For $K-1\rightarrow K$, by induction
    \begin{align*}
        &\norm{\exp{\i M_B(t_K)\Delta t}\cdots\exp{\i M_B(t_0)\Delta t}-
        \exp{\i\tilde{M}_B(t_K)\Delta t}\cdots\exp{\i\tilde{M}_B(t_0)\Delta t}
        }[\Lnorm]\\
        &\leq
        \norm{\left[\exp{\i M_B(t_K)\Delta t}-
        \exp{\i\tilde{M}_B(t_K)\Delta t}\right]
        \exp{\i M_Bt_{K-1}\Delta t}
        \cdots\exp{\i M_B(t_0)\Delta t}}[\Lnorm]\\
        &+
        \bigg\|\exp{\i\tilde{M}_B(t_K)\Delta t}
        \times\\
        &\left[\exp{\i M_B(t_{K-1})\Delta t}\cdots\exp{\i M_B(t_0)\Delta t}-
        \exp{\i\tilde{M}_B(t_{K-1})
        \Delta t}\cdots\exp{\i\tilde{M}_B(t_0)\Delta t}\right]\bigg\|_{\Lnorm}\\
        &\leq|\Delta t|K \norm{M_B(t)-\tilde{M}_B(t)}[\Lnorm],
    \end{align*}
    where by definition $|\Delta t|K=|t|$. Thus, we can estimate for the ordered
    exponential
    \begin{align*}
        &\norm{\adj{\oe{\i M_B}{0}{t}}-\adj{\oe{\i\tilde{M}_B}{0}{t}}}[\Lnorm]\\
        &\leq |t|\norm{M_B(t)-\tilde{M}_B(t)}[\Lnorm]
        \leq CJ^2t^2\exp{-a\{l/3-v|t|\}}.
    \end{align*}
    We define the local operator $B$ as
    \begin{align*}
        B:=\frac{1}{N\sqrt{2\pi q}}\int_{-\infty}^\infty
        \adj{\oe{\i\tilde{M}_B}{0}{t}}\exp{-\frac{t^2}{2q}}\d t.
    \end{align*}
    Estimating as in \eqref{eq:trunc} for $T=\frac{l}{6v}$, we obtain
    \begin{align*}
        &\norm{\OB-B}[\Lnorm]\leq
        \norm{\int_{|t|\leq T\ldots}+\int_{|t|>T}\ldots}[\Lnorm]
        \leq
        C(1/N)J^2\exp{-al/3}\\
        &\times\left(\exp{avT}\frac{1}{\sqrt{2\pi q}}\int_{|t|\leq T}
        t^2\exp{-\frac{t^2}{2q}}\d t+
        \frac{1}{\sqrt{2\pi q}}\int_{|t|>T}t^2\exp{av|t|-\frac{t^2}{2q}}\d t
        \right).
    \end{align*}
    In analogy to \eqref{eq:trunc}, for the first term we obtain
    \begin{align*}
        \frac{1}{\sqrt{2\pi q}}\int_{|t|\leq T}
        t^2\exp{-\frac{t^2}{2q}}\d t\leq \frac{q}{2}.
    \end{align*}
    For the second term, applying integration by parts, we obtain
    \begin{align*}
        &\frac{1}{\sqrt{2\pi q}}\int_{|t|>T}t^2\exp{av|t|-\frac{t^2}{2q}}\d t=
        \sqrt{\frac{2}{\pi q}}\int_{T}^\infty
        t^2\exp{avt-\frac{t^2}{2q}}\d t\\
        &=\sqrt{\frac{2}{\pi q}}\bigg(
        Tq\exp{avT-\frac{T^2}{2q}}+q\int_{T}^\infty \exp{avt-\frac{t^2}{2q}}
        dt
        +avq\int_{T}^\infty t\exp{avt-\frac{t^2}{2q}}\d t\bigg).
    \end{align*}
    And hence
    \begin{align*}
        &\frac{1}{\sqrt{2\pi q}}\int_{|t|>T}t^2\exp{av|t|-\frac{t^2}{2q}}\d t
        \leq\exp{avt-\frac{T^2}{2q}}\frac{\sqrt{2}}{\sqrt{\pi}|1-avq|}
        \sqrt{q}\left(\sqrt{\frac{q\pi}{2}+T}\right).
    \end{align*}
    The final estimate is thus
    \begin{align*}
        \norm{\OB-B}[\Lnorm]\leq
        C_3J^2\max\set{q^{1/2},\, q}\exp{-3c_3l},
    \end{align*}
    for appropriate constants $C_3>0$, $c_3>0$. This completes the proof.
\end{proofEnd}

We now have all the ingredients for the main result.
\begin{theoremEnd}[normal]{theorem}\label{thm:obolor}
    Under Assumption \ref{ass:op}, there exist local, bounded and
    self-adjoint operators $L=L(j,\,l)$,
    $B=B(j,l)$,
    $R=R(j,l)$ with norms bounded by 1,
    such that for some constants $C_4>0$, $c_4>0$
    independent of $d$, $j$ or $l$
    \begin{align}\label{eq:finallemma}
        \norm{P_0-BLR}[\Lnorm]\leq
        C_4J^2\exp{-c_4l}.
    \end{align}
    The respective supports are $\X_{1,j}$,
    $\X_{j-3l-2,j+3l+3}$ and $\X_{j+1,d}$. The operator $B$ can be chosen
    w.l.o.g.\ to be positive.
\end{theoremEnd}
\begin{proofEnd}
    The operators $L$ and $R$ were defined in Lemma \ref{lemma:ol} and
    their properties follow therefrom. The operator $B$ was defined in Lemma
    \ref{lemma:ob}. W.l.o.g.\ we can assume it is positive, otherwise
    the same arguments as in \cite[Lemma 4]{Hastings} apply.

    By Lemmas \ref{lemma:tpq}, \ref{lemma:ob} and since
    $\norm{LR}[\Lnorm]\leq 1$,
    we obtain an error bound with
    asymptotic dependence on $l$ of the form $l^2\exp{-cl}$. Hence, we
    can pick constants $C_4>0$, $c_4>0$ to satisfy \eqref{eq:finallemma}.
    This completes the proof.
\end{proofEnd}

\bibliographystyle{acm}
\bibliography{literature}

\appendix
\section{Proofs for \Cref{sec:alaw}}\label{app:proofs}
\printProofs

\end{document}